\documentclass [reqno]{amsart}
\usepackage {amsmath,amssymb}
\usepackage {color,graphicx, dsfont}
\usepackage [utf8]{inputenc}
\usepackage [dvips,centering,includehead,width=14cm,height=22.5cm]{geometry}

\newcommand {\ignore}[1]{}

\newcommand {\Mrs}{M_{(r,s)}}
\newcommand {\Nrs}{N_{(r,s)}}
\newcommand {\ori}{{\text {\rm or}}}

\DeclareMathOperator {\ev}{ev}
\DeclareMathOperator {\ft}{ft}
\DeclareMathOperator {\mult}{mult}

\sloppypar

\renewenvironment {enumerate}%
  {\rule{1mm}{0mm}\begin {oldenumerate}%
    \parskip1ex plus0.5ex \itemsep 0mm \parindent 0mm}%
  {\end {oldenumerate}}

\renewenvironment {itemize}%
  {\rule{1mm}{0mm}\begin {olditemize}%
    \parskip1ex plus0.5ex \itemsep 0mm \parindent 0mm}%
  {\end {olditemize}}

\newlength {\sidetextwidth}
\newlength {\sidepicwidth}
\newsavebox {\sidepicbox}
\newenvironment {sidepic}[1]{%
  \par%
  \def\message##1{}%
  \hbadness=10000%
  \setlength{\sidetextwidth}{\textwidth}%
  \savebox{\sidepicbox}{\scalebox{0.6}{\input {pics/#1}}}%
  \settowidth{\sidepicwidth}{\usebox{\sidepicbox}}%
  \addtolength{\sidetextwidth}{-6mm}%
  \addtolength{\sidetextwidth}{-\sidepicwidth}%
  \begin {minipage}{\sidetextwidth}}%
  {\end {minipage}%
  \hfuzz=10000pt%
  \hspace {4mm}%
  \begin {minipage}{\sidepicwidth}\usebox{\sidepicbox}\end {minipage}%
  \hfuzz=0pt%
  \hbadness=1000%
  }

\theoremstyle {plain}
\newtheorem {theorem}{Theorem}[section]

\newtheorem {lemma}[theorem]{Lemma}

\newtheorem {corollary}[theorem]{Corollary}
\theoremstyle {definition}
\newtheorem {definition}[theorem]{Definition}
\theoremstyle {remark}
\newtheorem {remark}[theorem]{Remark}
\newtheorem {example}[theorem]{Example}

\newtheorem {notation}[theorem]{Notation}

\newcommand{\NN}{\mathds{N}}
\newcommand{\ZZ}{\mathds{Z}}
\newcommand{\QQ}{\mathds{Q}}
\newcommand{\RR}{\mathds{R}}
\newcommand{\CC}{\mathds{C}}
\newcommand{\PP}{\mathds{P}}
\newcommand{\calM}{\mathcal{M}}
\newcommand{\calP}{\mathcal{P}}

\theoremstyle {plain}

\newenvironment{LG}{{\bf LG}.  \footnotesize}{}

\title {Refined broccoli invariants}
\author {Lothar G\"ottsche \and Franziska Schroeter}
\address {Lothar G\"ottsche, International Center For Theoretical Physics, Strada Costiera
11, 34151 Trieste, Italy}
  \email{gottsche@ictp.it}
\address {Franziska Schroeter, Fachbereich Mathematik (AD), Universität Hamburg, Bundesstr. 55, 20146 Hamburg, Germany}
  \email{franziska.schroeter@uni-hamburg.de}


\begin {document}

  \begin {abstract}
We introduce a tropical enumerative invariant depending on a variable $y$ which generalizes the tropical refined Severi degree. We show that this refined broccoli invariant is indeed independent of the point configuration, and that it specializes to a tropical descendant Gromov-Witten invariant for $y=1$ and to the corresponding broccoli invariant for $y=-1$. Furthermore, we define tropical refined descendant Gromov-Witten invariants which equal the corresponding refined broccoli invariants giving a new insight to the nature of broccoli invariants. We discuss various possible generalizations, e.g. to refinements of bridge curves and Welschinger curves.
  \end {abstract}

  \maketitle
  \section {Introduction} \label {sec-intro}

In \cite{BGo14} refined tropical enumerative invariants are introduced, which interpolate between Severi degrees and \textit{totally real Welschinger invariants}; the latter count real curves passing through only real points and not through pairs of complex conjugate points. For these so called \textit{tropical refined Severi degrees} 
$N_\text{trop}^{(X,L),\delta}(y)$  (Definition \ref{def-refinedS}) we consider certain tropical curves through a collection $\calP$ of points in tropical general position and count each 
curve $C$ with a refined (curve) multiplicity $\text{mult}(C;y)$, which is an expression in one variable $y$ (Definition \ref{def-refinedS}).  
The (refined) Severi degrees and the disconnected Welschinger invariants  considered here count possibly reducible curves (or more precisely maps from disconnected  curves), but in \cite{BGo14} in the same way 
irreducible tropical refined Severi degrees are introduced, which now interpolate  between the Gromov-Witten invariants and the connected totally real Welschinger invariants (both of which count maps from connected  curves).

The refined multiplicity $\text{mult}(C;y)$ has 
analogously to the usual Mikhalkin curve multiplicity $\text{mult}_\CC$ (Definition \cite[2.16]{Mik05}) the property that it is a product of vertex multiplicities $[\text{mult}_\CC (v)]_y$, 
which however now are not integers, but Laurent polynomials in $y^{1/2}$.
This allows to carry over the tropical proof of invariance of Gathmann and Markwig for tropical Gromov-Witten invariants \cite{GM05a} to show that the tropical refined Severi degree and the irreducible tropical refined Severi degree   
do not depend on the choice of points in $\calP$ as long as they remain in general position \cite{IM13}. 
In the rest of the paper we will differently from \cite{BGo14} mostly consider only irreducible tropical refined Severi degrees and just call them tropical refined Severi degrees, and by Welschinger invariants we mean connected Welschinger invariants.
The tropical refined Severi degree is a  symmetric Laurent polynomial in the variable $y$ and specializes to the corresponding tropical Gromov-Witten invariant $N_\text{trop}^{(X,L),\delta}$ for $y=1$ and to the corresponding tropical Welschinger invariant for $y=-1$. By the usual correspondence theorems e.g. \cite{Mik05}, it follows that these tropical invariants equal their counterparts in algebraic geometry, and hence that tropical refined Severi degrees interpolate between Gromov-Witten invariants and totally real Welschinger invariants \cite{BGo14}.

So far, these tropical refined Severi degrees are purely tropical objects and their nature is still unclear. 
In \cite{GoS12} \textit{refined invariants} $\widetilde N^{(S,L),\delta}(y)$ have been introduced, which can be viewed as a refinement of BPS invariants, and it was conjectured that for a $\delta$-very ample line bundle $L$ on a smooth toric surface $S$ they agree with the tropical refined Severi degrees. This has been shown for $\PP^2$ 
and rational ruled surfaces for small values of $\delta$ in \cite{BGo14} and \cite{GoK15}. In \cite{FS12} Filippini and Stoppa relate tropical refined Severi degrees to the wall-crossing formula of refined Donaldson-Thomas invariants.Furthermore, Mikhalkin \cite{Mik15} relates them to the weighted count of certain real curves. Finally, in \cite{NPS16} the authors give a geometric interpretation of the refined curve multiplicity using motivic measures of certain semialgebraic domains of Hilbert schemes of points.

It is unclear how the tropical refined Severi degrees are related to other known enumerative invariants, and if we can produce many new interesting tropical invariants when we specialize to other values of $y$.It might look mysterious at a first glance why they do interpolate only between Severi degrees and \textit{totally real} Welschinger invariants and not Welschinger invariants in general. At least one can hope to generalize their definition in order to obtain non-totally real tropical Welschinger invariants as specialization as well. By \cite{Shu06} the tropical curves $C$ we consider for non-totally real tropical Welschinger invariants, called \textit{Shustin curves}, are counted with multiplicities, which are \textit{not} products of vertex multiplicities. Shustin shows there their invariance by proving that they equal the corresponding non-totally real Welschinger invariants, but he could not carry over the tropical proof of Gathmann and Markwig.

This was so far only possible by the detour via \textit{broccoli curves} \cite{GMS13} which are certain rational tropical curves having again the property that their curve multiplicity $m_C$ is a product of vertex multiplicities $m_V$. The broccoli curves passing through a configuration $\calP$ of points in general position can be deformed in a prescribed way by the so called \textit{bridge algorithm} in order to obtain Shustin curves passing through $\calP$. Then the weighted number of broccoli curves through $\calP$ equals the weighted numbers of Shustin curves through $\calP$. Since it has been proven by the methods of Gathmann and Markwig that the weighted numbers of broccoli curves, the so called \textit{broccoli number} $N^B_{r,s}(\Delta,F,\calP)$ does not depend on the choice of $\calP$, it follows that the non-totally real tropical Welschinger invariant is also invariant. Also for broccoli invariants, their nature in algebraic geometry is unclear.

In this paper we will introduce \textit{refined broccoli curves} (Definition \ref{def-or}). Refined broccoli curves have a slightly simpler definition than broccoli curves \cite{GMS13}, the main difference being that there are no restrictions on the parity of the weights of the edges of a refined broccoli curve. They form a class of tropical curves which includes the set of tropical curves used for the definition of the tropical (refined) Severi degrees and the set of broccoli curves. 
Using the idea of proof of Gathmann and Markwig we can show that the count of refined broccoli curves yields an invariant (Theorem \ref{thm-invariance}), namely the \textit{refined broccoli invariant} $N^{rB}_{(r,s)}(y,\Delta,F)$ as defined in Definition \ref{def-or-refined-brocc-inv}. Refined broccoli invariants are again defined by counting the refined broccoli curves $C$ with a refined (curve) multiplicity $m_C(y)$, 
which is a product of vertex multiplicities $m_V(y)$. While the vertex multiplicities $m_V(y)$ are in general rational functions in $y^{1/2}$, the multiplicities $m_C(y)$ and thus the refined broccoli invariants turn out to be again symmetric Laurent polynomials in  $y$, which generalize the tropical refined Severi degrees of \cite{BGo14}. In addition to the mentioned interpolation between (tropical) Gromov-Witten invariants 
and (tropical) totally real Welschinger invariants, they interpolate for certain choices of $r$ and $s$ also between broccoli invariants for $y=-1$ (Corollary \ref{cor-minuseins}) and tropical descendant Gromov-Witten invariants for $y=1$ (Lemma \ref{cor-same} and Corollary \ref{dGW=rdGW}). This interpolation between tropical descendant Gromov-Witten invariants and (refined) broccoli curves uncovers a new relation of broccoli curves to already known tropical invariants.

Tropical descendant Gromov-Witten invariants $\tilde{N}_{\Delta,\mathbf{k}}^\text{trop}(\alpha)$ have been introduced in \cite{MR09} and \cite{BGM11} and can be considered as the tropical counterpart of descendant Gromov-Witten invariants (see e.g. \cite{FP97},\cite{KM98}). In \cite{BGM11} it has been proven that the tropical and non-tropical descendant Gromov-Witten invariants agree for $\PP^2$.  We also introduce \textit{refined descendant curves} (Definition \ref{def-trdi}) and show that for a given configuration of points $\calP$ in general position we have a bijection between refined broccoli curves through $\calP$ and refined descendant curves through $\calP$ (Lemma \ref{lem-equ}). This implies that the refined broccoli invariant equals the corresponding \textit{refined descendant Gromov-Witten invariant} (Corollary \ref{cor-same}), which is the weighted count of refined descendant curves.

In a forthcoming work \cite{GoSch} we will introduce floor diagrams for refined broccoli curves and prove a Caporaso Harris type recursion formula for the refined broccoli invariants. 
We will also give a formula for these invariants in terms of a Heisenberg algebra action.
We also make speculations about higher genus broccoli invariants, including a conjectural generating function.

\subsection*{Acknowledgments}
We would like to thank Eugenii Shustin and Peter Overholser for helpful discussions. The second author was partially supported during the research by GIF Grant no.\ 1174-197.6/2011, the Minkowski-Minerva Center for Geometry at the Tel Aviv University, by Grant no.\ 178/13 from the Israel Science Foundation, and by the RTG 1670 ``Mathematics Inspired by String Theory and Quantum Field Theory'' funded by the German Research Foundation (DFG). 

  \section {Preliminaries and reminder} \label {sec-pre}

\subsection {Non-oriented and oriented marked curves} \label {sub-marked}
For details of the following, the reader is referred to \cite[Section 2]{GMS13}. Let $r,s \geq 0$.
\subsubsection{Non-oriented marked curves and their moduli space}\label {subsub-nonorient} 
An $(r,s)$-marked (plane tropical rational parametrized) curve $C$ of degree $\Delta$ consists of an abstract, not necessarily
connected, metric graph $\Gamma$ and a continuous map $h: \Gamma \to \RR^2$, which is integer affine linear on edges of $\Gamma$, 
and such that the balancing condition is satisfied at vertices of $h(\Gamma)$. Furthermore, we require that each connected component of $\Gamma$ is rational, i.e.\ has first Betti number equal to zero. The collection of unbounded edges $(x_1,\ldots,x_{r+s})$ of $\Gamma$, that are mapped to $0$ by $h$ consists of the \textit{markings} of $C$ and the labeled set $(v(y_1),\ldots,v(y_n))$ of the direction vectors of the non-contracted unbounded edges (``\textit{ends of} $C$'') is called the \textit{degree} $\Delta$ of $C$. That is, our curves are \textit{labeled} as all unbounded edges, not only the markings, are labeled. $|\Delta|$ is the number of elements in $\Delta$. Note that if $\Delta$ contains respectively $d$ times the vector $\begin{pmatrix} -1,0 \end{pmatrix}$, $\begin{pmatrix} 0,-1 \end{pmatrix}$ and $\begin{pmatrix} 1,1 \end{pmatrix}$ we talk about \textit{degree} $d$ \textit{curves}. The \textit{combinatorial type} $\alpha$ of an $(r,s)$-marked curve $C$ is the data of $C$, but forgets about the length of the edges of $\Gamma$. Two marked curves $(\Gamma;h;x_1,\ldots,x_{r+s})$ and $(\Gamma';h';x'_1,\ldots,x'_{r+s})$ are \textit{isomorphic} if there is a homeomorphism $\varphi:\Gamma \to \Gamma'$, $x_i$ is mapped to $x'_i$ under $\varphi$, every edge of $\Gamma$ is mapped bijectively to an edge of $\Gamma'$ by an affine map of slope $\pm 1$, and we have $h' \circ \varphi =h$.\\
We say that an $(r,s)$-marked curve $C=(\Gamma,h)$ is \textit{connected} if $\Gamma$ is connected. 
The set all isomorphism classes of connected $(r,s)$-marked curves of degree $\Delta$ will be denoted by $\Mrs(\Delta)$ which  is the \textit{moduli space} $\calM^{\text{lab}}_{0,r+s,\text{trop}}(\RR^2,\Delta)$ of connected $(r+s)$-marked plane labeled rational tropical curves. This is a tropical variety, whose structure as abstract polyhedral complex is inherited from the combinatorial types of the curves in it. 
The open polyhedron in $\Mrs(\Delta)$ consisting of curves with the same combinatorial type $\alpha$ will be denoted by $\Mrs^\alpha$. It turns out that its dimension can be computed as $2$ plus the number of bounded edges in $h(\Gamma)$ of a curve of combinatorial type $\alpha$ \cite[Remark 2.5 and 2.7]{GMS13}.\\
As a convention we will draw an $(r,s)$-marked curve $C$ always as $h(\Gamma)\subset \RR^2$. An edge $e$ of $h(\Gamma)$ will be displayed in bold and called \textit{even} if the greatest common divisor of the entries of its direction vector 
$v(e)\in \ZZ^2$ (``\textit{weight of $e$}'') is $2$. Otherwise, the edge $e$ is called \textit{odd} and will be drawn as a thin line. When the curve has two ends of the same direction adjacent to the same vertex (``\textit{double ends}'') we will draw these edges parallel to each other with a small distance in between. The images of the markings $h(x_i)$ will be depicted as dots, which will be small if $i = 1,\ldots,r$ (``\textit{real markings}'') and big for $i=r+1,\ldots,r+s$ (``\textit{complex markings}'').\\

Note that our definitions and notations agree most of the time with \cite{BGo14} and \cite{BGM11}. One important difference is that in \cite{BGo14} the tropical curves are \textit{unlabeled} 
which means that only the markings are labeled and in \cite{BGM11} the tropical curves are {partially unlabeled}, that is, only the markings and some of the ``left ends'' are labeled 
\cite[Definition 3.1]{BGM11}. This has consequences for the enumeration of tropical curves. Namely, if $m$ unlabeled ends of the same direction vector are adjacent to the same 
vertex in a (partially) unlabeled curve, they are not distinguishable, therefore they contribute with an automorphism factor of $\frac{1}{m!}$ to the contribution of that curve in 
the count of (partially) unlabeled curves. When the same situation appears (for labeled ends) in a labeled curve the automorphism factor does not appear in the count of labeled 
curves. In most cases, we are interested in the count of unlabeled curves, which will however be performed by  counting labeled curves and then dividing by the number of different 
ways of labeling the ends, see for example Definition \ref{def-wbb}.

\subsubsection{Evaluation map and conditions in general position}\label {subsub-ev} 
Let $C$ be an $(r,s)$-marked curve of degree $\Delta=(v(y_1),\ldots,v(y_n))$ and let $F \subset \{1,\ldots,n\}$. The ends $y_i$ with $i \in F$ are called \textit{fixed ends} and they are indicated in drawings by a small orthogonal bar at the ``infinite 
side''; their role in enumerative problems will be illustrated in \cite{GoSch}. Furthermore, let $M\subset \Mrs(\Delta)$ be a polyhedral subcomplex. Then the \textit{evaluation map with respect to} $F$ \textit{and} $M$ is given by 

$$\ev_{F,M}:M \to (\RR^2)^{r+s} \times \prod_{i \in F} (\RR^2/\langle 
v(y_i)\rangle),~C \mapsto \big((h(x_1),\ldots,h(x_{r+s})),(h(y_i):~i\in F)\big).$$

In most cases $M= \Mrs(\Delta)$, so we can write $\ev_{F}$ instead of $\ev_{F,M}$. If $F=\emptyset$ then we write $\ev$ instead of $\ev_F$. $\ev_{F,M}$ is a morphism of polyhedral complexes, in particular it is continuous and linear on each polyhedron $\Mrs^\alpha (\Delta)$.\\ 
A \textit{collection of conditions for} $\ev_{F,M}$ is a tuple $\calP=\big((P_1,\ldots,P_{r+s}),(Q_i:~i\in F)\big)$ of points $P_i \in \RR^2$ and lines $Q_i \in \RR^2/\langle v(y_i) \rangle$. Hence $\ev_{F,M}^{-1}(\calP)$ is the set of $(r,s)$-marked curves in $M$ passing through $P_i$ at the marking $x_i$ and such that the fixed end $y_i$ is mapped to the line $Q_i$. The \textit{locus of conditions in general 
position} (for $\ev_{F,M}$) is the complement of the union $\bigcup_{\alpha} ev_{F,M}(M^\alpha \subset \Mrs^\alpha(\Delta))\subset \RR^{2(r+s)+|F|}$, where the union is taken over all combinatorial types $\alpha$ such that $\ev_{F,M}(M^\alpha)$ has dimension at most 
$2(r+s)+|F|-1$. As $\ev_{F,M}$ is linear on $\Mrs^\alpha(\Delta)$, this means that $\ev_{F,M}(M^\alpha)$ lies in the locus of conditions in general position if $M^\alpha$ has dimension at least $2(r+s)+|F|$. If $F =\emptyset$ we talk about \textit{points in general position}.
\begin{remark}[Difference in definition to \cite{GM05b}] \label{rem-olddef}
In contrast to \cite{GM05b} and also \cite{Mik05} marked curves passing through points in general position are here not supposed to be $3$-valent by definition.
A definition closer to \cite{GM05b} will be given in \ref{def-restricted-pts}. Note that we need this 
more general definition here since we are for instance interested in counts of Welschinger curves as defined in Definition \ref{def-wbb} (b).
\end{remark}
Some of the $(r,s)$-marked curves in top-dimensional cells of $\Mrs(\Delta)$ passing through points in general position have a property allowing to orient unmarked edges in a unique way. This property is a generalization of \cite[Remark 3.7]{GM05b} and is Lemma 2.13 (b) of \cite{GMS13}:
\begin{lemma} \label{lem-unique}
Let $M\subset \Mrs(\Delta)$ be a polyhedral subcomplex and let $\calP$ be a collection of conditions in general position. Assume there is a curve $C \in \ev_F^{-1}(\calP)$. If the combinatorial type of $C$ has dimension $2(r+s)+|F|$ and every vertex of $C$ that is not adjacent to a marking is $3$-valent, then every connected component of $\Gamma \setminus (x_1\cup\ldots \cup x_{r+s})$ has exactly one unmarked end $y_i$ with $i \notin F$. 
\end{lemma}
\begin{example}\label{ex-gegen_gen}
The condition about the valency of vertices is important; for instance in the example below we have $r=1$, $s=2$, $F=\emptyset$, and the points are in general position. However, there is one connected component of $\Gamma \setminus (x_1\cup x_2 \cup 
x_3)$ which has two unmarked ends. Although the points are in general position, the big dot lying on the edge contributes only with one dimension to the space of conditions, which is therefore in total $4$-dimensional. But we need actually $|\Delta|-1=5$ conditions to ensure that there is exactly one curve passing through the points. In our situation there is a $1$-parameter family of curves passing through them.
\begin {center} \begin{picture}(0,0)%
\includegraphics{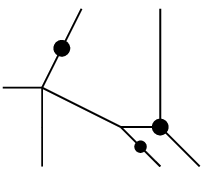}%
\end{picture}%
\setlength{\unitlength}{4144sp}%
\begingroup\makeatletter\ifx\SetFigFont\undefined%
\gdef\SetFigFont#1#2#3#4#5{%
  \reset@font\fontsize{#1}{#2pt}%
  \fontfamily{#3}\fontseries{#4}\fontshape{#5}%
  \selectfont}%
\fi\endgroup%
\begin{picture}(924,744)(4669,-3403)
\end{picture}%
 \end {center}
\end{example}

Lemma \ref{lem-unique} implies Remark 2.14 of \cite{GMS13} which we state here for its importance as 
\begin{corollary} \label{cor-orient}
Let $C \in \ev_F^{-1}(\calP)$ as in Lemma \ref{lem-unique}. Then there is a unique way to orient all unmarked edges of $C$ such that the orientation points towards the unique unmarked non-fixed end of the component of $\Gamma \setminus (x_1 \cup \ldots \cup x_{r+s})$ containing the edge.
\end{corollary}
This orientation of fixed ends will then be always inwards (the curve). The orientation of an edge will be indicated on pictures by an arrow.
\subsubsection{Oriented marked curves and their relationship to unoriented marked curves}\label {subsub-orient}

An \textit{oriented} $(r,s)$-marked curve of degree $\Delta$ consists of the data of an unoriented $(r,s)$-marked curve of degree $\Delta$ and a unique orientation of each unmarked edge of $C$.  The \textit{set of fixed ends} $F \subset \{1,\ldots,n\}$ \textit{of} $C$ contains all the indices $i$ of those unmarked ends $y_i$ of $C$ which are oriented inwards. $\Mrs^\ori(\Delta,F)$ denotes the set of all connected oriented $(r,s)$-marked curves of degree $\Delta$ and  set of fixed ends $F$ and we will drop $F$ in this notation if $F=\emptyset$. Note that there is a natural forgetful map $\ft:\Mrs^\ori(\Delta,F) \to \Mrs(\Delta)$ which forgets the orientation of edges and allows to compare both spaces, see Remark 2.16 \cite{GMS13}. In particular, $\ft$ is a morphism of polyhedral complexes and is injective on each cell of $\Mrs^\ori(\Delta,F)$.\\

In \cite{GMS13} the authors turned their attention to a particular class of curves in $\Mrs^\ori(\Delta,F)$, which they describe in terms of vertex types. This is because if an unoriented curve $C\in \Mrs(\Delta)$ satisfies the conditions of Lemma \ref{lem-unique} and if we can make it into an oriented curve $C'\in \Mrs^\ori(\Delta,F)$ composed only of vertex types (1) - (7) below such that the indices of the fixed ends 
of $C$ form the set $F$ of $C'$, then the orientation of $C'$ agrees with the natural orientation of $C$ given by Corollary \ref{cor-orient}, see \cite[Lemma 2.20]{GMS13}. Let us recall briefly these vertex types.

\begin {sidepic}{vertices-a} 
\begin{definition}[Old vertex types] \label{def-vertextypes} A \textit{vertex type} of an oriented $(r,s)$-marked curve $C$ is the information of the number, parity (even or odd) and the orientation of the adjacent edges of a vertex in $C$. In the following list, a vertex of type (6) is of type (6b) if the two odd edges are parallel, otherwise of type (6a). The two odd edges in type (6b) must be unmarked ends if they appear in a 
Welschinger curve as defined in definition \ref{def-wbb}. Also in case (8) the two odd edges must be unmarked ends. The two odd edges for (9) must not be parallel, which is indicated by the arc (as in (6a)).\end{definition}
\end{sidepic}
If $a$ denotes the (complex) Mikhalkin multiplicity \cite[2.16]{Mik05} of the vertex $V$ of a certain type of the list, then we assign a new vertex multiplicity $m_V$ to this vertex, which depends on the vertex type. In the figure, $\text{i}$ denotes the imaginary unit.  Note that as Mikhalkin's definition is only for $3$-valent vertices (of $h(\Gamma)$), $a$ is computed in case (8) as the absolute value of the determinant of the direction vectors of the even edges. Denote by $w(y_i)$ the weight of the end $y_k$ of the connected oriented marked curve $C$. Then its multiplicity $m_C$ is given by
\begin{equation}\label{eq:brocmult}
m_C=\prod_{k=1}^n \text{i}^{w(y_k)-1}\prod_{V \text{ in }C}m_V.
\end{equation}
If the curve $C$ is non-connected, then its multiplicity is the product of the multiplicities of the connected components of $C$.

\subsection {Broccoli, Welschinger and bridge curves}

The main achievement of \cite{GMS13} was the introduction of \textit{broccoli curves}, which can be seen as auxiliary tropical curves and allow to prove the invariance of Welschinger numbers \cite{Wel03, Wel05} for the count of real (rational) curves through real points and complex conjugate points in general position by purely tropical means, even if a proof via correspondence theorems was known before \cite{Shu06}. 
The community is interested in such a proof since for instance the tropical invariance proof \cite{IKS09} for Welschinger numbers associated to real curves through only real points in general position gave rise to tropical Welschinger invariants for real curves of higher genus and recursive formulas for these numbers. By \cite[Lemma 5.8]{GMS13} broccoli curves and also \textit{Welschinger curves} are special cases of a particular class of $(r,s)$-marked curves, which we call \textit{bridge curves}. Broccoli curves and Welschinger curves can be deformed into each other by a series of bridge curves lying on a so called \textit{bridge graph} \cite[Remark 5.12]{GMS13}. For broccoli and Welschinger 
curves there exists an oriented and an unoriented version of the definition, while for bridge curves there is only an oriented one.  By \cite[Proposition 3.3]{GMS13} and \cite[Proposition 4.10]{GMS13} there is a bijection between oriented broccoli/Welschinger curves through a configuration $\calP$ of conditions in general position and unoriented broccoli/Welschinger curves through $\calP$.

Here is the oriented version of the definitions. 

\begin {definition}[Welschinger, broccoli and bridge curves] \label {def-wbb}
Let $C \in \Mrs^\ori(\Delta,F)$ be an oriented $(r,s)$-marked curve. $C$ is called 
\begin{itemize}
\item[(a)] \textit{broccoli curve} \cite[Definition 3.1.(a)]{GMS13}, if it is only composed of vertices of type (1) - (6) from Definition \ref{def-vertextypes}.
\item[(b)] \textit{Welschinger curve} \cite[Definition 4.6.(a)]{GMS13}, if it is only composed of vertices of type (1) - (5), (6b), (7), or (8) from Definition \ref{def-vertextypes}.
\item[(c)] \textit{bridge curve}, if it consists of vertices of type (1) - (8) and at most one vertex of type (9), and there is a bijection between the vertices of type (7) and those of type (8) or (9) satisfying suitable conditions (specified in  \cite[Definition 5.2]{GMS13}).
\end{itemize}
\end {definition}

Let $\Mrs^B(\Delta,F)$ be the closure of broccoli curves in $\Mrs^\ori(\Delta,F)$, which is a polyhedral subcomplex of $\Mrs^\ori(\Delta,F)$. Assume now $F$ satisfies $r+2s+|F|=|\Delta|-1$ to ensure that $\Mrs^B(\Delta,F)$ is non-empty \cite[Lemma 2.21]{GMS13}. When $\calP \in \RR^{2(r+s)+|F|}$ is a collection of conditions in general position for $\ev_{F,\Mrs^B(\Delta,F)}: \Mrs^B(\Delta,F) \to 
\RR^{2(r+s)+|F|}$, then the \textit{broccoli invariant} (still with respect to \ $\calP$) is defined as 

\begin{equation}\label{eq:broccount} \Nrs^B(\Delta,F,\calP)=\frac{1}{|G(\Delta,F)|}\sum_{C} m_C,\end{equation}

where $G(\Delta,F)$ is the symmetric subgroup of $S_n$ of all permutations $ \sigma$ satisfying $\sigma(i)=i$ for all $i\in F$ and $v(y_i)=v(y_{\sigma(i)})$ for all $i=1,\ldots,n$; the sum is taken over all broccoli curves $C$ in $\Mrs^B(\Delta,F)$ such that $\ev_{F,\Mrs^B(\Delta,F)}(C)=\calP$. The symmetric group factor is necessary since we deal with labeled curves and so we overcount by $|G(\Delta,F)|$. By Theorem 3.6 of \cite{GMS13} $\Nrs^B(\Delta,F,\calP)$ does not depend on 
$\calP$:

\begin{theorem}[Invariance for broccoli] \label{old-thm-invariance}
The broccoli invariants $\Nrs^B(\Delta,F,\calP)$ are independent of the collection $\calP$ of conditions and therefore we write them from now on as $\Nrs^B(\Delta,F)$ or as $\Nrs^B(\Delta)$ if $F=\emptyset$.
\end{theorem}

The idea of proof will be briefly repeated in the beginning of section \ref{sec-inv} and follows the idea of proof of the main theorem of \cite{GM05a}.

Analogously, let $\Mrs^W(\Delta,F)$ be the closure of Welschinger curves in $\Mrs^\ori(\Delta,F)$ and assume that $F$ is chosen such that $r+2s+|F|=|\Delta|-1$ to make sure that the polyhedral subcomplex $\Mrs^W(\Delta,F)$ of $\Mrs^\ori(\Delta,F)$ is non-empty \cite[Lemma 2.21]{GMS13}. When $\calP \subset \RR^{2(r+s)+|F|}$ is a collection of conditions in general position for $\ev_{F,\Mrs^W(\Delta,F)}: \Mrs^W(\Delta,F) \to 
\RR^{2(r+s)+|F|}$, then the tropical \textit{Welschinger number} is defined as

$$\Nrs^W(\Delta,F,\calP)=\frac{1}{|G(\Delta,F)|}\sum_C m_C,$$

where the sum is taken over all Welschinger curves $C$ in $\Mrs^W(\Delta,F)$ with 
$\ev_{F,\Mrs^W(\Delta,F)}(C)=\calP$.  
In \cite[Corollary 5.17]{GMS13} it is shown that for $\Delta$ a \textit{toric del Pezzo degree} (see \cite[Definition~4.22]{GMS13} for the definition, note that in this case the associated toric surface $X(\Delta)$ is $\PP^2$, $\PP^1\times\PP^1$ or the blowup of $\PP^2$ in at most 3 general points) the tropical Welschinger numbers $\Nrs^W(\Delta,F,\calP)$ are independent of collection $\calP$ of conditions. Therefore we will in this case just write them as $\Nrs^W(\Delta,F)$, or as $\Nrs^W(\Delta)$ if $F=\emptyset$.

Note that this definition of Welschinger numbers agrees with Shustin's definition \cite{Shu06} if $F = \emptyset$ and $\Delta$ contains only primitive vectors, see \cite[Remark 4.19]{GMS13}. Nevertheless, the definition of a Welschinger curve \cite[Remark 4.14]{GMS13}, and of its multiplicity is different \cite[Lemma 4.18]{GMS13}, even if $F=\emptyset$ and all vectors in $\Delta$ are primitive. For example, Shustin's curves are unparametrized and are always nodal, see Definition \ref{def-tropical-curve}.(b). The multiplicity of a Welschinger 
curve $C$ in Shustin's sense is \textit{not} a product of vertex multiplicities, but depends on more global information of the curve, see \cite[Definition 4.15]{GMS13}. But in the case 
where we consider curves with no complex marking, that is $s=0$, Shustin's curve multiplicity agrees with Mikhalkin's \textit{real curve multiplicity} $\mult_\RR(C)$ \cite[Definition 7.19]
{Mik06} and is given as follows: 
$$\mult_\RR(C)=\begin{cases} 0, & \text{ if }\mult_\CC(C)\text{ is even},\\
				1, & \text{ if }\mult_\CC(C)\equiv 1 \text{ mod } 4  ,\\
				-1, & \text{ if }\mult_\CC(C)\equiv 3 \text{ mod } 4,
				\end{cases}$$
where $\mult_\CC(C)$ is the usual \textit{Mikhalkin multiplicity of the curve} $C$, i.e.\ the product of the Mikhalkin vertex multiplicities, taken over all $3$-valent vertices in $C$. 

\subsection {Tropical refined Severi degrees and their properties}
Let $\Delta$ be a lattice polytope, $X=X(\Delta)$ the associated projective toric surface and $L=L(\Delta)$ the associated tautological toric line bundle. Then the count of $\delta$-nodal curves in $|L|$ through $\dim |L|-\delta$ points in general position, which do not contain a toric boundary divisor as a component, equals the count of certain tropical curves dual to $\Delta$ through points in tropical general position, if we count every tropical curve with a multiplicity. The invariants we obtain in this way are called \textit{toric Severi degrees} $N^{(X,L),\delta}$ of $X$ and $L$. 
In \cite{BGo14} the authors define \textit{tropical refined Severi degrees} $N_{\text{trop}}^{(X,L),\delta}(y)$. By \cite[Theorem 4.3]{BGo14} this invariant $N_{\text{trop}}^{(X,L),\delta}(y)$ equals for suitable choices of $(X(\Delta), L(\Delta))$ the \textit{refined invariant} $\tilde{N}^{(X,L),\delta}(y)$ \cite[Definition 2.4]{BGo14}, i.e.\ it appears as coefficient of a certain generating function \cite[Equation (2.3)]{BGo14} introduced in \cite{GoS12}, which involves the $\chi_{-y}$-genus of certain relative Hilbert schemes.\\ 

Let us briefly review the definition of tropical refined Severi degrees. 
\begin{definition}[Simple/nodal tropical curves] \label{def-tropical-curve}
Let $C$ an $n$-marked connected (plane tropical rational parametrized) curve of degree $\Delta$ as defined in \cite[Definition 4.1]{GKM07}. 
\begin{itemize}
\item[(a)]$C$ is called \textit{simple} \cite[Definition 4.2]{Mik05} if the map $h: \Gamma \to \RR^2$ is an immersion, all the vertices of $\Gamma$ are $3$-valent, for any $x \in \RR^2$ the preimage $h^{-1}(x)$ consists of at most $2$ points and if $a, b \in \Gamma$ with $a 
\neq b$ such that $h(a) = h(b)$, then $a$ and $b$ must not be a vertex of $\Gamma$.
\item[(b)] $C$ is called \textit{nodal} \cite[Definition 3.4.(3)]{BGo14} if its dual Newton subdivision consists only of triangles and parallelograms. The number of nodes $\delta$ of $C$ is given by $|\mathring{\Delta} \cap \ZZ^2| - g(C)$, where $\mathring{\Delta}$ is the interior of $\Delta$ and the \textit{genus} $g(C)$ of $C$ is defined as the first Betti number of the underlying graph $\Gamma$ .   
\end{itemize}
\end{definition}
\begin{remark}\label{rem-Delta}
The letter $\Delta$ will be used simultaneously for lattice polytopes and degrees of tropical curves. There is no risk of confusion,  since any tropical degree $\Delta$ defines a lattice polytope when we rotate all vectors $v_i$ in $\Delta$ by $-\pi/2$ and draw them in 
$\ZZ^2$ one after the other starting at a lattice point of $\ZZ^2$, each with lattice length equal to the weight of the corresponding edge, in a chain. Producing a tropical degree from a lattice polytope works too if we assume that the tropical curve of that degree has only unmarked ends of weight $1$.
\end{remark}
\begin{definition}[Tropical refined Severi degree]\label{def-refinedS}
Let $\delta \geq 0$, $\Delta$  a lattice polytope, $X=X(\Delta)$ the associated projective toric surface and $L=L(\Delta)$ the associated line bundle as before. Let $C$ be a simple $(|\Delta \cap \ZZ^2|-1-\delta)=(\dim |L|-\delta)$-marked $\delta$-nodal curve of degree $\Delta$. Note that here, differently from above for broccoli curves, we consider unlabeled curves. 
Furthermore, let $\text{mult}_\CC(V)=a\in \NN$ be the Mikhalkin multiplicity \cite[Definition 2.16]{Mik05} of a $3$-valent vertex $V$ of $C$ and $\calP$ a configuration of $|\Delta \cap \ZZ^2|-1-\delta$ points in $\RR^2$ in (tropical) general position.
\begin{itemize}
\item[(a)]  The \textit{refined vertex multiplicity} \cite[Definition 3.5]{BGo14} is defined as $$[\text{mult}_\CC(V)]_y=\frac{y^{a/2}-y^{-a/2}}{y^{1/2}-y^{-1/2}}=y^{ (a-1)/2 } + y^{ (a-3)/2} +\ldots+y^{ -(a-1)/2}.$$
\item[(b)] The \textit{refined curve multiplicity} \cite[Definition 3.5]{BGo14} of $C$ is given by $$\text{mult}(C;y)=\prod_{V \text{ is } 3\text{-valent in }C} [\text{mult}_\CC(V)]_y.$$
\item[(c)] The \textit{tropical refined Severi degree} \cite[Definition 3.7]{BGo14} associated to $X$ and $L$ with respect to $\calP$ is defined as $$N_{\text{trop}}^{(X,L),\delta}(y,\calP)=\sum_{C}\text{mult}(C;y),$$ where the sum is taken over all simple $(|\Delta \cap \ZZ^2|-1-\delta)$-marked $\delta$-nodal curves $C$ of degree $\Delta$ passing through $\calP$. 
Note again that these are unlabelled curves. Because of this, differently from \eqref{eq:broccount}, this formula does not contain a prefactor accounting for the permutations of the labelings. 
\end{itemize}
\end{definition}
\begin{remark} \label{rem-severi_ori}
A curve $C$ considered for a tropical refined Severi degree passes through points in general position, and so by \cite[Remark 3.7]{GM05b} there is a unique way to orient the edges of $C$ such that arrows on the edges always show away from the marked points and in the direction of the only possible end.
\end{remark}

By \cite[Theorem 1]{IM13} $N_{\text{trop}}^{(X,L),\delta}(y,\calP)$ does not depend on $\calP$, so we can drop $\calP$ in the notation. Furthermore, $N_{\text{trop}}^{(X,L),\delta}(y)$ is a Laurent polynomial in $\NN[y,y^{-1}]$, which is invariant under the action $y \rightarrow y^{-1}$. More importantly, the tropical refined Severi degree interpolates between the Severi degree and the corresponding Welschinger number in the following sense:

\begin{theorem}[Interpolation]\label{thm-interpolate}
Let $\delta,X,L,\Delta$ and $C$ be as in Definition \ref{def-refinedS}. Then we have $\text{mult}(C;1)=\text{mult}_\CC(C)$ and $\text{mult}(C;-1)=\text{mult}_\RR(C)$. Furthermore, we have $N_{\text{trop}}^{(X,L),\delta}(1)=N^{(X,L),\delta}$.\\ 
Remember that we denote the interior of $\Delta$ by $\mathring{\Delta}$ and denote by $\partial \Delta$ the boundary of $\Delta$. Considering tropical curves of genus $0$: if $X(\Delta)$ is equal to $\PP^2$, $\PP^1 \times \PP^1$ or $\PP^2$ blown up in up to three real points equipped with the standard real structure we have $N_{\text{trop}}^{(X,L),|\mathring{\Delta} \cap \ZZ^2 |}(-1)=W^{(X,L)}$, where $W^{(X,L)}$ is the Welschinger number associated to $X(\Delta)$ and $L(\Delta)$ counting real rational curves in $X(\Delta)$ through $|\partial \Delta \cap \ZZ^2 |-1$ real points in general position.
\end{theorem}
\begin{remark} \label{rem-alsocomplexpoints}
Note that tropical refined Severi degrees $N_{\text{trop}}^{(X,L),\delta}(y)$ do not specialize to Welschinger invariants counting curves which also pass through configurations of points containing pairs of complex conjugate points. 
\end{remark}

  \section {Refined broccoli curves and their count} \label {sec-broc}
\ignore{In the following we are not considering arbitrary toric surfaces, but only those 
associated to $h$-tranverse polygons $\Delta$. This is due to the fact that we can count 
curves on these surfaces via \textit{floor diagrams} \cite{BM08}, which are for example 
used to prove the existence of \textit{node polynomials} for Severi degrees in the sense 
of Ardila and Block \cite{AB13}.
\begin{definition}[$h$-transverse polygon] \label{def-htransverse}
A lattice polygon $\Delta$ is called $h$-\textit{transverse} \cite[section 2.2]{BM08} if 
the direction vectors of its edges are all of the shape $(\alpha,0)$ (``vertical edge''), 
$(0,\alpha)$ (``horizontal edge'') or $(\alpha,\pm 1)$ with $\alpha \in \ZZ$. 
\end{definition}
\begin{example}
***Special cases of $h$-transverse polygons are the standard $d$-simplex ($->$ define 
tropical curves of degree $d$), Hirzebruch and so on.***
\end{example}
}
\subsection{Oriented (rational) refined broccoli curves} \label{sub-obroc}
We start with the definition of one of the main objects in this paper which is a refined version of broccoli curves. 
\begin{definition}[Oriented refined broccoli curves] \label{def-or}
An oriented $(r,s)$-marked curve $C\in \Mrs^\ori(\Delta,F)$ is called \textit{oriented refined broccoli curve} if it is composed of vertices of the following three vertex types, where any parity of the adjacent edges is allowed (as long as the balancing condition is satisfied). The corresponding \textit{refined vertex multiplicity} is written below each vertex. Here, $a$ denotes again the Mikhalkin multiplicity $\text{mult}_\CC(V)$ of the vertex $V$. Note that for the vertex type (I) the marking is real and the special case (6b) in the list in Definition \ref{def-vertextypes} of type (III) has multiplicity $\frac{2}{y^{1/2}+y^{-1/2}}$.
\begin {center} \scalebox{0.7}{\input {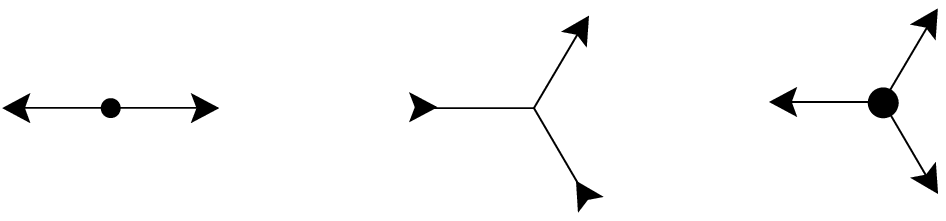}} \end {center}
Let $y_1,\ldots, y_n$ be the ends of $C$.
For a fixed end $y_{i}$ of $C$ (i.e. $i\in F$), let $m_{y_i}(y):=\frac{y^{w(y_i)/2}+(-1)^{w(y_i)}y^{-w(y_i)/2}}{y^{1/2}+(-1)^{w(y_i)}y^{-1/2}}$, and for a non-fixed end $y_i$, let $m_{y_i}(y):=\frac{y^{w(y_i)/2}-(-1)^{w(y_i)}y^{-w(y_i)/2}}{w(y_i)(y^{1/2}-(-1)^{w(y_i)}y^{-1/2})}$.
We define the \textit{refined multiplicity} of $C$ as  
$$m_C(y)=\prod_{i=1}^n m_{y_i}(y)\prod_{V \in C} m_V(y).$$
\end{definition}
\begin{remark}
The vertex type (II) is the same as considered in Definition \ref{def-refinedS} (see also Remark \ref{rem-severi_ori} for orientations), whereas type (III) is completely new. Type (I) was not considered as a vertex in \cite{BGo14}. However, for consistency reasons with 
\cite{GMS13} we continue to use it. It is not clear from the definition that the refined curve multiplicity $m_C(y)$ is a Laurent polynomial 
with similar properties as the tropical refined Severi degree we have defined in Definition \ref{def-refinedS}, but we will prove it in Lemma \ref{lem-welldefined} below. Note that the refined vertex multiplicity $m_V(y)$ is in general not a Laurent polynomial. 
\end{remark}

\begin{definition}[Broccoli index of a refined oriented broccoli curve]\label{def-brocindex}
Let  $C\in \Mrs^\ori(\Delta,F)$ be an oriented refined broccoli curve. Let $V_{cm}$ be the set of vertices of $C 
\subset \RR^2$ of even Mikhalkin multiplicity and such that a complex marking is adjacent 
to each of them. Analogously, let $V_{wcm}$ be the set of vertices of $C$ of even 
Mikhalkin multiplicity without complex marking. Let $E_{f}$ be the set of even fixed ends of $C$, i.e. $E_{f}=\big\{i\in F\bigm| w(y_i)\hbox{ even} \big\}$, and let
$E_{n}$ be the non-fixed even ends of $C$. We also write $e_f=\#E_f$ and $e_n=\#E_n$. Then we define the \textit{broccoli index of $C$} as 
$$i_B(C)= -\# V_{cm}-e_{f}+\#V_{wcm}+e_{n}.$$
\end{definition}
\begin{remark}\label{rem-even}
Note that the Mikhalkin multiplicity of a vertex is even if and only if  at least one of the adjacent edges is even. 
\end{remark}
\begin{notation}[Old broccoli curve]
We call a broccoli curve in the sense of \cite{GMS13} an \textit{old broccoli curve}.
\end{notation}
\begin{lemma}[Broccoli index for an old broccoli curve]\label{lem-brocindex}
Let $C$ be an old broccoli curve. Then we have the following equality for the number $n_{(x)}$ of vertices of type $(x)$ in a broccoli curve $C$ in the sense of \cite{GMS13}:
$$n_{(3)}+n_{(4)}+e_n=n_{(6)}+e_{f}.$$
Hence, $i_B(C)=0$ if $C$ is an old broccoli curve.
\end{lemma}
\begin{remark}
It follows from Lemma \ref{lem-welldefined} below that the converse is also true: if $i_B(C)=0$ for an oriented refined broccoli curve $C$, then $C$ is an old broccoli curve. This direction is of particular importance for enumeration of curves: if there was a curve $C$ with $i_B(C)=0$, which is not an old broccoli curve, then Lemma \ref{lem-welldefined} would imply $m_C(-1) \neq 0$. 
\end{remark}
\begin{proof}[Proof of Lemma 3.6.]
The argument is similar to the one of the proof of Lemma 4.18 \cite{GMS13}. 

Let $\widetilde{C}$ be a connected component (\textit{``broccoli''}) of the subgraph $C_{\text{even}}\subset C$ consisting of even edges. $\widetilde{C}$ 
is a $3$-valent graph after forgetting about markings and considering $1$-valent vertices in the boundary of $\widetilde{C}$ as ends. Let $\widetilde e_f$ and $\widetilde e_n$ respectively be the number of even fixed ends and even non-fixed ends of $C$ which are ends of $\widetilde C$. The number of ends of $\widetilde C$ is equal to $\widetilde{n}_{(3)}+\widetilde{n}_{(6)}+\widetilde e_n+\widetilde e_f $, so there are $\widetilde{n}_{(3)}+\widetilde{n}_{(6)}+\widetilde e_n+\widetilde e_f -2$ vertices in $\widetilde{C}$. On the other hand, this number of vertices in $\widetilde{C}$ is also 
equal to $\widetilde{n}_{(4)}$. As by \cite[3.1.(b)]{GMS13} $\widetilde{C}$ contains exactly one stem, and the stems are precisely the vertices of type $(3)$ and the even non-fixed ends of $C$, we get $\widetilde{n}_{(3)}+\widetilde e_n=1$.
This gives the equation $\widetilde{n}_{(6)}+\widetilde{e}_f=\widetilde{n}_{(4)}+1=\widetilde{n}_{(4)}+\widetilde{n}_{(3)}+\widetilde e_n$ for the number of vertex types and ends of $\widetilde{C}$. Taking into account all connected components of $C_{\text{even}}$ we obtain the result.
\end{proof}
\begin{example}[Counterexample]\label{ex-counter}
Not every refined oriented broccoli curve is a broccoli curve in the old sense. For example, the following three marked curves through points in general position are composed of the vertices of type (I) and (II), but they contain forbidden vertices as depicted in the proof of \cite[Proposition 3.3]{GMS13}, so they are not broccoli curves (problematic vertices are marked by the orientation of their adjacent edges).
\begin {center} \scalebox{1.3}{\begin{picture}(0,0)%
\includegraphics{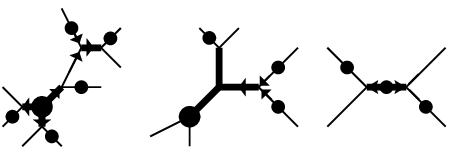}%
\end{picture}%
\setlength{\unitlength}{4144sp}%
\begingroup\makeatletter\ifx\SetFigFont\undefined%
\gdef\SetFigFont#1#2#3#4#5{%
  \reset@font\fontsize{#1}{#2pt}%
  \fontfamily{#3}\fontseries{#4}\fontshape{#5}%
  \selectfont}%
\fi\endgroup%
\begin{picture}(2049,654)(4489,-2953)
\end{picture}%
} \end {center}
Comparing the list of vertex types allowed for a broccoli curve in Definition \ref{def-wbb} with the list in Definition \ref{def-or} we see that these examples actually contain already all problematic vertex types, namely:
\begin {center} \scalebox{1.3}{\input {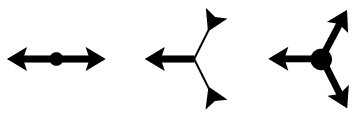}} \end {center}
\end{example}

The class of refined broccoli curves is larger than the class of broccoli curves considered 
in \cite{GMS13}, but their relationship is clarified by the following  
\begin{corollary}[Relationship with old broccoli curves] \label{cor-old}

\begin{itemize}
\item[(a)] If $C$ is an old  broccoli curve, then it is an oriented refined broccoli curve. In this case its multiplicity $m_C$ equals the specialization $m_C(-1)$ of its refined multiplicity at $y=-1$.
\item[(b)] Conversely every oriented refined broccoli curve $C$ that is not an old broccoli curve satisfies $m_C(-1)=0$. 
\end{itemize}
\end{corollary}
\begin{proof}
\begin{itemize}
\item[(a)] As the class of allowed vertex types is larger for oriented refined broccoli curves than for old broccoli curves the first statement is true. Similar to the expansion of $\frac{y^{a/2}-y^{-a/2}}{y^{1/2}-y^{-1/2}}=y^{(a-1)/2} +y^{(a-3)/2}+\ldots+y^{ -(a-1)/2}$ of Definition \ref{def-refinedS}, we have 
\begin{equation} \label{otherex}
\frac{y^{a/2}-(-1)^a y^{-a/2}}{y^{1/2}+y^{-1/2}}= y^{ (a-1)/2} -y^{(a-3)/2}+\ldots +(-1)^{a-1} y^{ -(a-1)/2}.\end{equation}
With this we see that the multiplicity of a vertex of type (II) gives the right specialization at $y=-1$ if $a$ is odd. The same is true for a vertex of type (III) if $a$ is odd. So it remains to discuss vertices of type (3), (4) and (6), and the contribution of fixed and non-fixed ends. By Lemma 
\ref{lem-brocindex} we have the equation $n_{(3)}+n_{(4)}+e_n=n_{(6)}+e_f$ for an old broccoli curve.  We define $N:=n_{(6)}+e_f$. 
The contribution of the vertices of type (3), (4) and (6) and the ends in $E_n$ and $E_f$ is 
\begin{align*}
&\textstyle \prod\limits_{j=1}^{N-e_n} \frac{y^{a_j/2}-y^{-a_j/2}}{y^{1/2}-y^{-1/2}}\prod\limits_{j\in E_n}\frac{y^{w(y_j)/2}-y^{-w(y_j)/2}}{w(y_j)(y^{1/2}-y^{-1/2})}
\prod\limits_{i=j}^{N-e_f}\frac{y^{\widetilde{a}_j/2}+y^{-\widetilde{a}_j/2}}{y^{1/2}+y^{-1/2}}\prod\limits_{j\in E_f}\frac{y^{w(y_j)/2}+y^{-w(y_j)/2}}{y^{1/2}+y^{-1/2}}\\& \textstyle
=\prod\limits_{j=1}^{N-e_n} \frac{y^{a_j/2}-y^{-a_j/2}}{y^{1/2}+y^{-1/2}}\prod\limits_{j\in E_n}\frac{y^{w(y_j)/2}-y^{-w(y_j)/2}}{w(y_j)(y^{1/2}+y^{-1/2})}
\prod\limits_{j=1}^{N-e_f}\frac{y^{\widetilde{a}_j/2}+y^{-\widetilde{a}_j/2}}{y^{1/2}-y^{-1/2}}\prod\limits_{j\in E_f}\frac{y^{w(y_j)/2}+y^{-w(y_j)/2}}{y^{1/2}-y^{-1/2}},
\end{align*}
where $a_j$ is the (even) Mikhalkin multiplicity of a vertex of type (3) or (4) and $\widetilde{a}_j$ is the (even) Mikhalkin multiplicity of a vertex of type (6). 

We now see that factors in the first, second, third and forth products give  $a_j \text{i}^{a_j-1}$, $\text{i}^{w(y_j)-1}$, $\text{i}^{\widetilde{a}_j-1}$ and $\text{i}^{w(y_j)-1}$ respectively, if we plug in $y=-1$. Here we use that 
$\frac{y^{\widetilde{a}_j/2}+y^{-\widetilde{a}_j/2}}{y^{1/2}-y^{-1/2}}=\frac{y^{1/2}(y^{\widetilde{a}_j } +1) } { y^{\widetilde{a}_j/2}(y-1)}$ for $y \notin \{0,1\}$ and the  expansion formula \eqref{otherex} above.\\
Denoting $O_{n}$ and $O_f$ respectively the indices of the non-fixed odd edges and the fixed odd edges of $C$, the multiplicity $m_C(y)$ is obtained from the above product by multiplying by 
$$\prod_{j\in O_n}\frac{y^{w(y_j)/2}+y^{-w(y_j)/2}}{w(y_j)(y^{1/2}+y^{-1/2})}\prod_{j\in O_f}\frac{y^{w(y_j)/2}-y^{-w(y_j)/2}}{y^{1/2}-y^{-1/2}}.$$
Which specializes at $y=-1$ to $\prod_{j\in O_f\cup O_n}\text{i}^{w(y_j)-1}.$
Comparing with the Definition \ref{def-vertextypes} of the broccoli multiplicity $m_C$, gives the claim.
\item[(b)] In Lemma \ref{lem-welldefined} (b) below we will show that the refined curve multiplicity $m_C(y)$ of such a curve $C$ is in $\ZZ[y,y^{-1}]$ and has a factor $y+2+y^{-1}$. Thus plugging in $y=-1$ yields $0$.
\end{itemize}
\end{proof}

\begin{remark}\label{rem-ver}
Note that we do not obtain the right specialization for old broccoli curves locally, i.e.\ for every vertex separately, but we have to use the equation of Lemma \ref{lem-brocindex}.  
\end{remark}

\begin{remark}\label{other}
In Definition \ref{def-or} we could also have made a different choice for the multiplicities associated to fixed and non-fixed ends.
A simpler choice would be to put $m'_{y_j}(y)=\frac{y^{w(y_j)/2}-y^{-w(y_j)/2}}{w(y_j)(y^{1/2}-y^{-1/2})}$ for $y_j\not \in F$ and 
$m'_{y_j}(y)=\frac{y^{w(y_j)/2}+y^{-w(y_j)/2}}{y^{1/2}+y^{-1/2}}$ for $y_j \in F$. Denote $m'(C)$ the multiplicity of the curve $C$ obtained in this way. With the notations above, a modification of the proof of Corollary \ref{cor-old}(a) shows that for an old broccoli curve $C$ we get  that $m'_C(-1)=m_C\frac{\prod_{j\in  O_f} w(y_j)}{\prod_{j\in O_n}w(y_j)}$.
\end{remark}

\begin{lemma}[Properties of the curve multiplicity] \label{lem-welldefined} Let $C$ be an oriented refined broccoli curve and let $\overline{N}$ be the set of indices of the non-fixed ends of $C$.
\begin{itemize}
\item[(a)]  The refined curve multiplicity $m_C(y)$ is a symmetric Laurent polynomial in $y$, more precisely $m_C(y)\prod_{i\in \overline{N}} w(y_i)\in \ZZ[y^{\pm 1}]$ and $m_C(y)=m_C(y^{-1})$. In particular, if $\overline{N}=\emptyset$, then $m_C(y)\in \ZZ[y^{\pm 1}]$.
\item[(b)] If $C$ is an old broccoli curve, then $m_C(-1)\ne 0$.
\item[(c)] If $C$ is not an old broccoli curve, then $m_C(-1)=0$. More precisely  in this case $i_B(C)$ is even and strictly positive, and we can write $m_C(y)=(y^{1/2}+y^{-1/2})^{i_B(C)}f(y)$ for $f\in\QQ[y^{\pm 1}]$ with $f(-1)\ne 0$.
\end{itemize}
\end{lemma}
\begin{proof}
(0) Note that  $m_C(y)$ is a product of factors of the forms  $\frac{y^{a/2}-y^{-a/2}}{y^{1/2}-y^{-1/2}}$, $\frac{y^{a/2}+y^{-a/2}}{y^{1/2}+y^{-1/2}}$, for positive integers $a$, and each of these factors $f(y)$ satisfies the symmetry $f(y)=f(y^{-1})$. We will show  $m_C(y)\prod_{i\in \overline{N}} w(y_i)\in \ZZ[y^{\pm 1}]$. Then it follows automatically that $m_C(y)$ is a symmetric Laurent polynomial.

(1) First we assume that $C$ is an old broccoli curve and show (a) and (b) in this case.
When $V$ is a vertex of the (generally non-connected) subgraph $C_{\text{odd}}$ consisting of odd edges, then its Mikhalkin multiplicity $a$ is odd. If $V$ is either of type (II) or type (III) we can expand its refined vertex multiplicity $m_V(y)$ to obtain a Laurent 
polynomial in $\ZZ[y^{\pm 1}]$. In the same way, if $y_i$ is an end of $C$ belonging to $C_{odd}$, then its weight is odd and therefore 
$m_{y_i}(y)w(y_i)\in  \ZZ[y^{\pm 1}]$, and  $m_{y_i}(y)\in \ZZ[y^{\pm 1}]$ if $i\not \in \overline{N}$.

Let us now consider  $C_{\text{even}}$. By Lemma \ref{lem-brocindex} we have again the equality $n_{(3)}+n_{(4)}+e_n=n_{(6)}+e_f$ for the number of vertices of each vertex type and the number of even ends. Define again $N:=n_{(6)}+e_f$. 
We denote $(a_1,\ldots, a_N)$ the collection of the  Mikhalkin multiplicities of the vertices of type (3) or (4) and of the weights of the even non-fixed ends of $C$, and we denote $(\widetilde{a}_1,\ldots, \widetilde{a}_N)$ the collection of the Mikhalkin multiplicities of vertices of type (6) and of the weights of even fixed ends of $C$.
So we have to show that
$$\prod_{i=1}^N \frac{y^{a_i/2}-y^{-a_i/2}}{y^{1/2}-y^{-1/2}}\prod_{i=1}^N 
\frac{y^{\widetilde{a}_i/2}+y^{-\widetilde{a}_i/2}}{y^{1/2}+y^{-1/2}} \in \ZZ[y,y^{-1}].$$
We compute
$$\textstyle \prod\limits_{i=1}^N \left(\frac{y^{a_i/2}-y^{-a_i/2}}{y^{1/2}-y^{-1/2}}\frac{y^{\widetilde{a}_i/2}+y^{-\widetilde{a}_i/2}}{y^{1/2}+y^{-1/2}} \right)
=\prod\limits_{i=1}^N \left(\frac{y^{(a_i+\widetilde{a}_i)/2}-y^{-(a_i+\widetilde{a}_i)/2}}{y-y^{-1}}+\frac{y^{(a_i-\widetilde{a}_i)/2}-y^{-(a_i-\widetilde{a}_i)/2}}{y-y^{-1}}\right).$$
Clearly for an integer $n$ we have $\frac{y^n-y^{-n}}{y-y^{-1}}=\pm(y^{n-1}+y^{n-2}+\ldots +y^{-n+1}) \in \ZZ[y,y^{-1}]$.
As all multiplicities $a_i$ and $\widetilde{a}_i$ are even, this implies that  both $\frac{y^{(a_i+\widetilde{a}_i)/2}-y^{-(a_i+\widetilde{a}_i)/2}}{y-y^{-1}}$ and $\frac{y^{(a_i-\widetilde{a}_i)/2}-y^{-(a_i-\widetilde{a}_i)/2}}{y-y^{-1}}$
lie in $\ZZ[y,y^{-1}].$ This shows (a) for old broccoli curves. 

(2) Part (b) is straightforward from the definitions. Let again $y_1,\ldots, y_n$ be the ends of $C$. By Definition \ref{def-vertextypes}, we see that all the vertex multiplicities $m_V$ of an old broccoli curve $C$ are nonzero, and clearly the weights $w(y_i)$ of all the ends of $C$ are nonzero.
Thus the multiplicity of the curve is $m_C=\prod_{i=1}^n{\textit i}^{w(y_i)-1}\prod_{V\in C} m_V\ne 0$. On the other hand by 
Corollary  \ref{cor-old}, we have that $m_C(-1)= m_C\ne 0.$

(3) Now we want to show part (a) and also part (c) for general oriented refined broccoli curves $C$.
 The main idea here is to turn a refined oriented broccoli curve $C$, which is not an old broccoli curve, into an old broccoli curve $\widetilde{C}$ in a prescribed way (\textit{``broccolization'' of the curve}) and then to compare their broccoli indices.\\
Suppose $C$ is a refined oriented broccoli curve with one or several vertices of type (a), (b) or (c) of Example \ref{ex-counter}. Replace these vertices by vertices allowed for old broccoli curves following the algorithm below. Note that the curve we obtain in this way is not necessarily connected anymore.
\begin {itemize}
 \item[(i)] Replace a vertex (a) by two vertices of type (6) as depicted below.
 What is meant here (and  the cases (ii) and (iii) below will be similar) is the following.
 We choose a small disk $D$ in $\RR^2$ with center at the vertex (a). Inside this disk we replace the vertex (a) by the two vertices of type (6), so that the even edges connect as before to the rest of the curve. The four new odd edges are unbounded edges. If $v$ is the direction vector of the even edge of a new  vertex of type (6), we choose the two unbounded direction vectors $v_1$, $v_2$ to form a decomposition $-v=v_1+v_2$ of $-v$ into primitive vectors. This is possible, because we can use an integral linear transformation of $\RR^2$ to reduce ourselves to the case
 the $v=(n,0)$ for some positive integer $n$, and then we can choose e.g. $v_1=(-1,n)$, $v_2=(-n+1,-n)$. Note that (as also in cases (ii) and (iii) below) the ends of the new curve 
 are the same as the ends of $C$ plus some (in this case four) non-fixed ends of weight 1.
There will be two more vertices in $V_{cm}$, so we have $i_B(C')=i_B(C)-2 \#(a)$, where $C'$ 
is the curve we obtain by replacing all vertices (a) in this way.
 
\begin {center} \scalebox{1.3}{\begin{picture}(0,0)%
\includegraphics{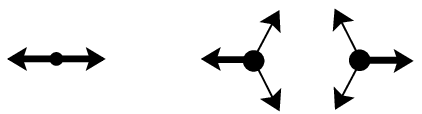}%
\end{picture}%
\setlength{\unitlength}{4144sp}%
\begingroup\makeatletter\ifx\SetFigFont\undefined%
\gdef\SetFigFont#1#2#3#4#5{%
  \reset@font\fontsize{#1}{#2pt}%
  \fontfamily{#3}\fontseries{#4}\fontshape{#5}%
  \selectfont}%
\fi\endgroup%
\begin{picture}(1903,493)(4468,-3642)
\put(5063,-3424){\makebox(0,0)[lb]{\smash{{\SetFigFont{8}{9.6}{\familydefault}{\mddefault}{\updefault}{\color[rgb]{0,0,0}$\longrightarrow$}%
}}}}
\end{picture}%
} \end {center}
  \item[(ii)] Replace a vertex (b) by one vertex of type (6) as depicted below. 
Here we mean that we choose a small disk $D$ in $\RR^2$ with center at the vertex (b). Inside this disk we replace the vertex (b) by a vertex of type (6), so that the edges have the same direction vectors. The even edge connects as before to the rest of the curve. The odd edges adjacent to the vertex of type (6) are unbounded. The odd edges that were before connected to the vertex (b) become unbounded edges. There 
will be one vertex less in $V_{wcm}$ and one more in $V_{cm}$, hence $i_B(C')=i_B(C)-2 
\#(b)$, where $C'$ is the curve we obtain by replacing all vertices (b) in this way.

  \begin {center} \scalebox{1.3}{\begin{picture}(0,0)%
\includegraphics{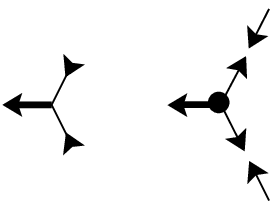}%
\end{picture}%
\setlength{\unitlength}{4144sp}%
\begingroup\makeatletter\ifx\SetFigFont\undefined%
\gdef\SetFigFont#1#2#3#4#5{%
  \reset@font\fontsize{#1}{#2pt}%
  \fontfamily{#3}\fontseries{#4}\fontshape{#5}%
  \selectfont}%
\fi\endgroup%
\begin{picture}(1243,899)(5920,-3843)
\put(6391,-3436){\makebox(0,0)[lb]{\smash{{\SetFigFont{8}{9.6}{\familydefault}{\mddefault}{\updefault}{\color[rgb]{0,0,0}$\longrightarrow$}%
}}}}
\end{picture}%
} \end {center}
  \item[(iii)] Replace a vertex (c) by three vertices of type (6) as depicted below, as follows.
 We choose a small disk $D$ in $\RR^2$ with center at the vertex (c). Inside this disk we replace the vertex (c) by  three vertices of type (6), so that the even edges connect as before to the rest of the curve. The six new odd edges are unbounded edges. Again, if $v$ is the direction vector of the even edge of a new  vertex of type (6), we choose the two unbounded direction vectors $v_1$, $v_2$ to form a decomposition $-v=v_1+v_2$ of $-v$ into primitive vectors. 
There will be two vertices more  in $V_{cm}$, hence 
$i_B(C')=i_B(C)-2 \#(c)$, where $C'$ is the curve we obtain by replacing all vertices (c) 
in this way.
\begin {center} \scalebox{1.3}{\begin{picture}(0,0)%
\includegraphics{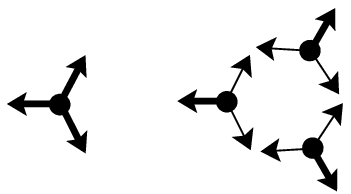}%
\end{picture}%
\setlength{\unitlength}{4144sp}%
\begingroup\makeatletter\ifx\SetFigFont\undefined%
\gdef\SetFigFont#1#2#3#4#5{%
  \reset@font\fontsize{#1}{#2pt}%
  \fontfamily{#3}\fontseries{#4}\fontshape{#5}%
  \selectfont}%
\fi\endgroup%
\begin{picture}(1601,861)(4904,-3815)
\put(5401,-3436){\makebox(0,0)[lb]{\smash{{\SetFigFont{8}{9.6}{\familydefault}{\mddefault}{\updefault}{\color[rgb]{0,0,0}$\longrightarrow$}%
}}}}
\end{picture}%
} \end {center}
\end{itemize}
In total, we see that the transformation of $C$ into an old broccoli curve $\widetilde{C}$ 
decreases the broccoli index, i.e.\
$$i_B(\widetilde{C})=i_B(C)-2\cdot \#((a)+(b)+(c)).$$
But since $\widetilde{C}$ is an old broccoli curve, it follows $i_B(C)> 0$ and $i_B(C)$ is 
even. Furthermore, since $\widetilde{C}$ is an old broccoli curve, its 
refined multiplicity lies in $\ZZ[y,y^{-1}]$  by part (a).\\

{\it Claim:} We can write $m_C(y)\prod_{i\in \overline{N}}w(y_i)=(y+2+y^{-1})^{i_B(C)/2}f(y)$ for $f \in \ZZ[y,y^{-1}]$ with $f(-1)\ne 0$.

By the above we see that by applying $i_B(C)/2$ times   (i), (ii) or (iii) to $C$ above we obtain an old broccoli curve $\widetilde{C}$, and by definition $\widetilde{C}$ satisfies the claim (with $i_B(\widetilde{C})=0$). Conversely we can view $C$ as having been obtained from $\widetilde{C}$  by applying $i_B(C)/2$ times the inverse of (i), (ii) and (iii).
Note that the steps (i), (ii), (iii) and their inverses do not change the ends of the curve except for possibly adding or removing some non-fixed ends of weight $1$, which have no influence on the multiplicity of the curve. 
For simplicity we denote 
$$\widetilde m_C(y):=m_C(y)\prod_{i\in \overline{N}} w(y_i).$$

We will use several times the following fact:
Consider $a_i,b_i,c_i\in \ZZ$ with $b_i$ odd and $c_i$ even. If the function
$$g(y)=\prod_{i=1}^{N_1}\frac{y^{a_i/2}-y^{-a_i/2}}{y^{1/2}-y^{-1/2}} \prod_{i=1}^{N_2}\frac{y^{b_i/2}+y^{-b_i/2}}{y^{1/2}+y^{-1/2}}\prod_{i=1}^{N_3}\frac{y^{c_i/2}+y^{-c_i/2}}{y^{1/2}+y^{-1/2}}$$
lies in  $\ZZ[y,y^{-1}]$, then it is divisible by $\prod_{i=1}^{N_3}(y^{c_i/2}+y^{-c_i/2})$ in $ \ZZ[y,y^{-1}]$.
This is true because for $c_i$ even, the factors $(y^{c_i/2}+y^{-c_i/2})$ are not divisible by $y^{1/2}+y^{-1/2}$ or $y^{1/2}-y^{-1/2}$.

(i) Let $C'$ be obtained from $C$ by applying (i), let $a$ be the (even) Mikhalkin multiplicity (which we can choose to  be equal) of the two new vertices of $C'$, and assume $\widetilde m_{C'}(y)\in \ZZ[y,y^{-1}]$. By the above observation, we can write $\widetilde m_{C'}(y)=(y^{a/2}+y^{-a/2})^2m_0(y)$ with $m_0(y)\in \ZZ[y,y^{-1}]$.
Then 
$$\widetilde m_C(y)=\left(\frac{y^{1/2}+y^{-1/2}}{y^{a/2}+y^{-a/2}}\right)^2\widetilde m_{C'}(y)=(y^{1/2}+y^{-1/2})^2 m_0(y)\in \ZZ[y,y^{-1}].$$

(ii) Let $C'$ be obtained from $C$ by applying (ii), let $a$ be the (even) Mikhalkin multiplicity of the changed vertex and assume $\widetilde m_{C'}(y)\in \ZZ[y,y^{-1}]$.
Then we can write $\widetilde m_{C'}(y)=(y^{a/2}+y^{-a/2}) m_0(y)$ with $m_0(y)\in \ZZ[y,y^{-1}]$, and we have
\begin{align*}
\widetilde m_C(y)&=\frac{y^{a/2}-y^{-a/2}}{y^{1/2}-y^{-1/2}}\frac{y^{1/2}+y^{-1/2}}{y^{a/2}+y^{-a/2}}\widetilde m_{C'}(y)=(y^{1/2}+y^{-1/2})^2\frac{y^{a/2}-y^{-a/2}}{y-y^{-1}}m_0(y)
\end{align*}
and $\frac{y^{a/2}-y^{-a/2}}{y-y^{-1}}\in \ZZ[y,y^{-1}]$ because $a$ is even. Furthermore by $\frac{y^{b}-y^{-b}}{y-y^{-1}}=y^{b-1}+y^{b-3}+\ldots+y^{-b+3}+y^{-b+1}$, we get
$\frac{y^{a/2}-y^{-a/2}}{y-y^{-1}}\big|_{y=-1}=(-1)^{a/2-1} \frac{a}{2}\ne 0$.

(iii) Let $C'$ be obtained  from $C$ by applying (iii), let $a$ be the (even) Mikhalkin multiplicity of the old vertex and $a_1,a_2,a_3$ the (even) Mikhalkin multiplicities of the new vertices, and assume $\widetilde m_{C'}(y)\in \ZZ[y,y^{-1}]$.
Then we can write $$\widetilde m_{C'}(y)=m_0(y) \prod_{i=1}^3(y^{a_i/2}+y^{-a_i/2})$$ with $m_0(y)\in \ZZ[y,y^{-1}]$, and 
\begin{align*}
\widetilde m_C(y)&=\left(\prod_{i=1}^3\frac{y^{1/2}+y^{1/2}}{y^{a_i/2}+y^{-a_i/2}}\right)\frac{y^{a/2}+y^{-a/2}}{y^{1/2}+y^{-1/2}} \widetilde m_{C'}(y)\\
&=(y^{a/2}+y^{-a/2})(y^{1/2}+y^{-1/2})^2m_0(y)\in \ZZ[y,y^{-1}].
\end{align*}

Altogether, we see that the inverses of steps (i), (ii) and (iii) transform curves $C'$ with $\widetilde m_{C'}(y)\in \ZZ[y,y^{-1}]$ into curves $C$ with  $\widetilde m_{C}(y)\in \ZZ[y,y^{-1}]$, and furthermore if we can write $\widetilde m_{C'}(y)=(y^{1/2}+y^{-1/2})^k f_0(y)$ with $f_0(-1)\ne 0$, then we can write $\widetilde m_{C'}(y)=(y^{1/2}+y^{-1/2})^{k+2} f_1(y)$ with $f_1(-1)\ne 0$.
As $\widetilde{C}$ is obtained from our original $C$ by applying $i_B(C)/2$ times the inverses of (i), (ii) or (iii), the claim follows.
\end{proof}
\ignore{\begin{remark}
Note that this statement is true for curves passing through points in general position as 
well as through points in special position. 
\end{remark}}
\begin{definition}[Oriented refined broccoli invariant]\label{def-or-refined-brocc-inv}
Let $\Mrs^{rB}(\Delta,F)$ be the closure of oriented refined broccoli curves in $\Mrs^\ori(\Delta,F)$ and 
assume we have $r+2s+|F|=|\Delta|-1$ in order to ensure that $\Mrs^{rB}(\Delta,F)$ is 
non-empty \cite[Lemma 2.21]{GMS13}. If $\calP$ is a collection of points in general 
position for $\ev_{\Mrs^{rB}(\Delta,F)}: \Mrs^{rB}(\Delta,F) \to \RR^{2(r+s)+|F|}$, then the 
\textit{oriented refined broccoli invariant} is defined as
$$\Nrs^{rB}(y,\Delta,F,\calP)=\frac{1}{|G(\Delta,F)|}\sum_{C}m_C(y),$$
where the sum taken over all oriented refined broccoli curves $C$ in $\Mrs^{rB}(\Delta,F)$ 
with $\ev_{\Mrs^{rB}(\Delta,F)}(C)=\calP$ and $|G(\Delta,F)|$ was defined in Definition \ref{def-wbb}. 
\end{definition}
From Corollary \ref{cor-old} follows now directly the following

\begin{corollary}\label{cor-minuseins}
We obtain the broccoli invariant $\Nrs^B(\Delta,F,\calP)$ of \cite{GMS13} if we set $y=-1$ in the oriented refined broccoli invariant $\Nrs^{rB}(y,\Delta,F,\calP)$.
\end{corollary}

It will be proven in Theorem \ref{thm-invariance} that $\Nrs^{rB}(y,\Delta,F,\calP)$ does not depend on $\calP$ and is therefore an invariant. 

\subsection{Unoriented (rational) refined broccoli curves} \label {sub-nbroc}
By \cite[Proposition 3.3]{GMS13} oriented and unoriented broccoli curves through points in general position are in bijection. Once we have introduced unoriented refined broccoli curves we are interested in a similar statement for refined broccoli curves.\\

\begin{definition}[Unoriented refined broccoli curves] \label{def-nonor}
An $(r,s)$-marked curve $C\in \Mrs(\Delta)$ is called \textit{unoriented refined broccoli curve} if it is composed of the following three vertex types, where any orientation and parity of the adjacent edges is allowed (as long as the balancing condition is satisfied). Here $a$ denotes again the Mikhalkin multiplicity $\text{mult}_\CC(V)$ of the vertex $V$. Note that for the vertex type (I') the marking is real.
\begin {center} \scalebox{0.7}{\input {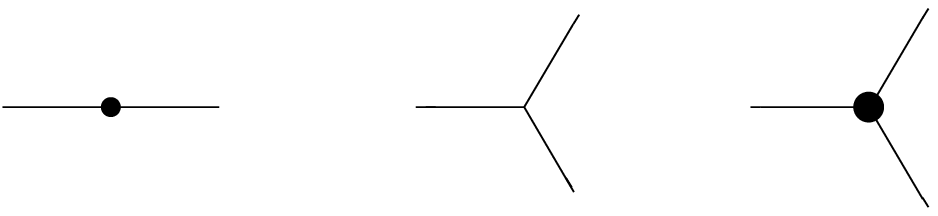}} \end {center}
We want to define multiplicities for non-oriented refined broccoli curves $C$, analogously to Definition \ref{def-or} for oriented ones. 
For this we need to specify which of the ends of $C$ are fixed.
Let $F\subset \{1,\ldots, n\}$ be a subset, such that $r+2s+|F|=|\Delta|-1$. 
For $i\in \{1,\ldots,n\}$ we put 
$$m_{y_i}(y):=\begin{cases} \frac{y^{w(y_i)/2}+(-1)^{w(y_i)}y^{-w(y_i)/2}}{y^{1/2}+(-1)^{w(y_i)}y^{-1/2}} & i\in F,\\
\frac{y^{w(y_i)/2}-(-1)^{w(y_i)}y^{-w(y_i)/2}}{w(y_i)(y^{1/2}-(-1)^{w(y_i)}y^{-1/2})} & i\not\in F.
\end{cases}$$
The multiplicity of $C$ with respect to $F$ is 
$$m_{C,F}(y):=\prod_{V\in C}m_V(y) \cdot \prod_{i=1}^n m_{y_i}(y).$$
Note that if all ends $C$ are simple, i.e.  $w(y_i)=1$ for all $i$, then $m_{C,F}(y)=\prod_{V\in C}m_V(y)$ is independent of $F$.
\end{definition}

\begin{lemma}[Relationship to oriented refined broccoli curves] \label{lemma-unori}
Let $r,s\geq 0$ and let $\Delta=(v_1,\ldots,v_n)$ be a collection of vectors in $\ZZ^2\setminus \{0\}$, and let $F\subset \{1,\ldots,n\}$ such that $r+2s+|F|=|\Delta|-1$. Furthermore, let $\calP=(P_1,\ldots,P_{2(r+s)+|F|})$ be a configuration of points in general position for $\ev_{F}$.Then there is a bijection between unoriented and oriented refined broccoli curves in $\Mrs(\Delta)$ respectively $\Mrs^\ori(\Delta,F)$ through $\calP$ with degree $\Delta$ and set of fixed ends $F$. 
\end{lemma}
\begin{proof}
We have to adapt the proof of Proposition 3.3 in \cite{GMS13}. Consider the forgetful map 
$\ft:\Mrs^\ori(\Delta,F) \to \Mrs(\Delta)$. Comparing vertex types (I), (II), (III) and 
(I'), (II'), (III'), it is clear that an oriented refined broccoli curve is mapped to 
an unoriented refined broccoli curve.
\begin{itemize}
 \item First, let us show that $\ft$ is injective on the set of 
curves through $\calP$ and let therefore $C \in \Mrs^\ori(\Delta,F)$ be an oriented 
refined broccoli curve through $\calP$. Since by \cite[Lemma 2.21]{GMS13} and the list of 
allowed vertex types the conditions for Lemma \ref{lem-unique} are satisfied, it 
follows by Corollary \ref{cor-orient} that $C$ has a unique natural orientation as 
described there, which is also the natural orientation of the unoriented 
refined broccoli curve $C'$ in the image of $C$ under $\ft$. Adapting the formulation and 
the proof of Lemma 2.20 \cite{GMS13} to our allowed vertex types, it follows that the 
orientation of $C$ has to agree with its natural orientation, hence two oriented refined 
broccoli curves cannot be mapped to the same unoriented refined broccoli curve under $\ft$.
\item $\ft$ is also surjective on the set of curves through $\calP$. Therefore, let $C' 
\in \Mrs(\Delta)$ be an unoriented refined broccoli curve. Since it is composed of 
vertices of types (I'), (II') and (III'), every vertex which is not adjacent to a complex 
marking is $3$-valent. By \cite[Proposition 2.11]{GM05b} the dimension of a combinatorial 
type can be computed as $|\Delta|-1+r-\sum_V(\text{val~}V-3)$, where the sum is taken 
over all vertices $V$ that are not adjacent to a complex marking. In our situation this 
dimension is equal to $2(r+s)+|F|$ and hence we can apply again Lemma \ref{lem-unique} to 
see that $C$ has a unique natural orientation. The vertices appearing in $C$ equipped 
with this natural orientation are then of type (I)-(III). 
\end{itemize}
\end{proof}

\begin{definition}[Unoriented refined broccoli invariant] 
\label{def-unor-refined-brocc-inv}
Let $F\subset \{1,\ldots, n\}$ be a subset, such that $r+2s+|F|=|\Delta|-1$.
Let $\Mrs^{urB}(\Delta)$ be the closure of unoriented refined broccoli curves in $\Mrs(\Delta)$. 
If $\calP$ is a collection of points in general 
position for $\ev_{F}$, then the 
\textit{unoriented refined broccoli invariant} with respect to $F$ is defined as
$$\Nrs^{urB}(y,\Delta,F,\calP)=\frac{1}{|G(\Delta,F)|}\sum_{C}m_C(y),$$
where the sum taken is over all unoriented refined broccoli curves $C$ in 
$\Mrs^{urB}(\Delta)$ with $\ev_{F}(C)=\calP$.
\end{definition}

\begin{corollary}[Equality of oriented and unoriented refined broccoli invariants]
Under the assumptions of Lemma  \ref{lemma-unori} we have
$\Nrs^{urB}(y,\Delta,F,\calP)=\Nrs^{rB}(y,\Delta,F,\calP).$
\end{corollary}

\subsection{Tropical refined descendant invariants} \label{sub-new} 

Descendant Gromov-Witten invariants can be defined for any genus on any smooth projective variety, making use of the virtual fundamental class 
on the moduli space of stable maps. 
Here we briefly recall the definition of rational descendant Gromov-Witten invariants of the complex projective plane.
For a nonnegative integer $n$ denote $M_{0,n}(\PP^2,d)$ the moduli space of $n$-marked rational stable maps of degree $d$ to $\PP^2$.
Its points correspond to tuples $(C,x_1,\ldots,x_n,f)$ of $C$ a connected complete rational curve with at most nodes as singularities, $x_1,\ldots,x_n$ distinct smooth points of $C$, $f:C\to \PP^2$ a morphism with $f_*[C]=dH$ where $H \in H_2(\PP^2,\ZZ)$ is the class of a line in $\PP^2$, such that the tuple $(C,x_1,\ldots,x_n,f)$ has only finitely many automorphisms.
For $i=1,\ldots,n$ we have  the evaluation morphisms $\ev_i:M_{0,n}(\PP^2,d)\to \PP^2$, sending $(C,x_1,\ldots,x_n,f)$ to $f(x_i)$.
Finally let $\psi_i\in A^1(M_{0,n}(\PP^2,d))$ be the first Chern class of the cotangent line bundle, i.e. the line bundle whose fiber over 
$(C,x_1,\ldots,x_n,f)$ is the cotangent space to $C$ at $x_i$. We denote $pt\in A^2(\PP^2)$ the class of a point.
For $a_1,\ldots,a_n\ge 0$ with $\sum_{i=1}^n (a_i+1)=3d-1$, the corresponding descendant Gromov-Witten invariant is defined as the intersection number
$$\langle \tau^{a_1}(pt) \cdots \tau^{a_n}(pt)\rangle_{d,0}:=\deg(\ev_1(pt)\psi_1^{a_1}\cdots \ev_n(pt)\psi_n^{a_n},[M_{0,n}(\PP^2,d)])\in \QQ.$$
(Note that we only consider stationary descendant Gromov-Witten invariants.) 
 
We now look at tropical descendant invariants. Although the original definition of tropical descendant Gromov-Witten invariants for rational curves goes back to \cite{MR09}, we stick here to the combinatorial definition of \cite{BGM11}, which agrees with the more intersection theoretical definition of \cite{MR09} in most cases, see Remark \ref{rem-old_descendant}.
\begin{definition}[Rational tropical descendant Gromov-Witten invariants] 
\label{def-dGW} 
Let $C$ be a labeled $m$-marked rational tropical curve of degree $\Delta=(v_1,\ldots,v_1,v_2,\ldots,v_2,\ldots,v_u\ldots,v_u)$, where the $v_i$ are distinct vectors and such that all vertices that are not adjacent to any of the markings are $3$-valent. Note that here again, we understand by ``labeled'' that all the ends and markings are labeled. The number of the $v_i$ in $\Delta$ does not have to be the same for all $i \in \{1,\ldots,u\}$. Assume $|\Delta|=n$. Choose a set $F\subset \{1,\ldots,n \}$ and define a sequence $\alpha=(\alpha_1,\alpha_2,\ldots)$ by $\alpha_i=\#\{v(y_j)\in \Delta|~j\in F \text{ and } v(y_j) \text{ has weight } i\}$.
We write $I^\alpha:=\prod_i i^{\alpha_i}$.

\begin{itemize}
\item The \textit{multiplicity} $m_C^{desc}$ of $C$ is defined as $\prod_{V} m_V$, where the product is taken over all $3$-valent vertices $V$ of $C$, which are not adjacent to a marking.
\item Let $\mathbf{k}=(k_0,k_1,\ldots)$ be a sequence of non-negative integers with only finitely many non-zero entries satisfying $m=|\mathbf{k}|:=\sum_i k_i$. 
Assume $\mathbf{k}$ satisfies $I(\alpha+\mathbf{k}):=0k_0+1(\alpha_1+k_1)+\ldots =|\Delta|-1-m$. Choose a vector 
$\mathbf{a}=(a_1,\ldots,a_{m})$ containing the number $i\in \NN$ exactly $k_i$ times -- in any order. Furthermore, let $\calP=(P_1,\ldots,P_{m})$ be a configuration of points in general position as defined in \cite[Definition 3.2]{MR09}. By \cite[Lemma 3.6]{MR09} this implies that every vertex of $\Gamma$ that is not $3$-valent has exactly one marking adjacent to it. Furthermore, we fix lines $Q_j\in \RR^2/\langle v(y_j)\rangle$ for $j \in F$ as in Subsection \ref{subsub-ev} such that the tuple $((P_1,\ldots,P_{m}),(Q_j|~j \in F))$ is a collection of conditions in general position. We define two types of descendant invariants:
\begin{itemize}
\item[(a)] $\widetilde{N}_{\Delta,\mathbf{k}}^\text{trop}(\alpha):=\frac{1}{I^\alpha}\frac{1}{|G(\Delta,F)|}\sum_C m_C^{desc} $, where the sum is taken over all $m$-marked rational tropical curves of degree $\Delta$ with markings $x_1,\ldots,x_m$ satisfying $h(x_i)=P_i$ for all $i=1,\ldots,m$, the marking $x_i$ is adjacent to a vertex in $\Gamma$ of valence $a_i+3$ for all $i=1,\ldots,m$, and the end $y_j$ for $j \in F$ is mapped to the line $Q_j\in \RR^2/\langle v(y_j)\rangle$, and $G(\Delta,F)$ is defined in Definition \ref{def-wbb}.
\item[(b)] $N_{\Delta,\mathbf{k}}^\text{trop}(\alpha):=\frac{1}{I^\alpha}\frac{1}{|G(\Delta,F)|}\sum_C m_C^{desc}$, where this time the sum goes over all $m$-marked rational tropical curves of degree $\Delta$ and with markings $x_1,\ldots,x_{m}$ such that $h(x_i)=P_i$ for $i=1,\ldots,m$, for every $i$ there are $k_i$ markings whose adjacent vertex has valence $i+3$, and the end $y_j$ for $j \in F$ is mapped to the line $Q_j\in \RR^2/\langle v(y_j)\rangle$, and $G(\Delta,F)$ is defined in Definition \ref{def-wbb}.
\end{itemize}
\end{itemize}
\end{definition}
\begin{remark}
\begin{itemize}
\item[(a)] Note that the considered curves are allowed to have marked vertices of valence bigger or equal to $4$.
\item[(b)] The numbers $\widetilde{N}_{\Delta,\mathbf{k}}^\text{trop}(\alpha)$ do not depend on the chosen vector $\mathbf{a}$.
\item[(c)] If $F=\emptyset$, then $\alpha$ is the zero sequence and we talk about \textit{absolute} rational tropical descendant Gromov-Witten invariants $\widetilde{N}_{\Delta,\mathbf{k}}^\text{trop}$ / $N_{\Delta,\mathbf{k}}^\text{trop}$. Otherwise we talk about \textit{relative} rational tropical descendant Gromov-Witten invariants. 
\item[(d)] In the  definition of $N_{\Delta,\mathbf{k}}^\text{trop}(\alpha)$ the order $j$ of a marking $x_j$ is not anymore related to the 
valency of the vertex to which $x_j$ is adjacent. That is, many more curves are considered in the second version. These two numbers are related by $N_{\Delta,\mathbf{k}}^\text{trop}(\alpha)=\frac{|\mathbf{k}|!}{\mathbf{k}!}\widetilde{N}^\text{trop}_{\Delta,\mathbf{k}}(\alpha)
$.
\end{itemize}
\end{remark}
\begin{remark}\label{rem-old_descendant}
\begin{itemize} 
\item[(a)] If $\Delta=(\underbrace{(-1,0),{\scriptstyle \ldots},(-1,0)}_{\alpha_1+\beta_1},\underbrace{(-2,0),{\scriptstyle \ldots} ,(-2,0)}_{\alpha_2+\beta_2},{\scriptstyle \ldots},\underbrace{(0,-1),{\scriptstyle \ldots},(0,-1)}_{d},\underbrace{(1,1),{\scriptstyle \ldots},(1,1)}_{d} ) $ for sequences $\alpha$, $\beta$ satisfying $I(\alpha+\beta)=d$ for some $d \in \NN_{>1}$, we write $\widetilde{N}_{d,\mathbf{k}}^\text{trop}(\alpha,\beta)$ / $N_{d,\mathbf{k}}^\text{trop}(\alpha,\beta)$ instead of $\widetilde{N}_{\Delta,\mathbf{k}}^\text{trop}(\alpha)$ / $N_{\Delta,\mathbf{k}}^\text{trop}(\alpha)$. Note that in this case, the vector $\mathbf{k}$ of Definition \ref{def-dGW} satisfies $I(\alpha+\beta+\mathbf{k})=3d-1+|\beta|-|\mathbf{k}|$. 
\item[(b)] The definition of $\widetilde{N}_{d,\mathbf{k}}^\text{trop}(\alpha,\beta)$ / $N_{d,\mathbf{k}}^\text{trop}(\alpha,\beta)$ in a) does not agree with Definition 3.7 \cite{BGM11}. Indeed we have $\widetilde{N}_{d,\mathbf{k}}^\text{trop}(\alpha,\beta)=\frac{1}{\beta!}\widetilde{N}_{d,\mathbf{k}}^\text{trop,BGM11}(\alpha,\beta)$ and $N_{d,\mathbf{k}}^\text{trop}(\alpha,\beta)=N_{d,\mathbf{k}}^\text{trop,BGM11}(\alpha,\beta)$. This is because in \cite{BGM11} the authors use partially labeled curves, that is for $\widetilde{N}_{d,\mathbf{k}}^\text{trop,BGM11}(\alpha,\beta)$ all left ends are labeled wheras for $N_{d,\mathbf{k}}^\text{trop,BGM11}(\alpha,\beta)$ only ends corresponding to the $\alpha$ sequence are labeled.
\item[(c)] If the vector $\mathbf{a}$ contains each $i\geq 0$ exactly $k_i$ times, then $\widetilde{N}_{d,\mathbf{k}}^\text{trop}$ equals the absolute descendant Gromov-Witten invariant defined in \cite{MR09} by \cite[Remark 3.3]{MR09}. This result has been generalized to relative descendant Gromov-Witten invariants in \cite{BGM11}.
\end{itemize}
\end{remark}
\begin{notation}
With ${\mathbf a}=(a_1,\ldots,a_m)$ as in Definition \ref{def-dGW}, we will also write 
$$\langle \tau_1^{a_1}(pt)\cdots\tau_m^{a_m}(pt)\rangle^{\text{trop}}_{\Delta,\alpha,0}:=\widetilde{N}_{\Delta,\mathbf{k}}^\text{trop}(\alpha).$$
For $F=\emptyset$, and therefore $\alpha$ the zero vector and $$\Delta=( \underbrace{(-1,0),\ldots,(-1,0)}_{d},\underbrace{(0,-1),\ldots,(0,-1)}_{d},\underbrace{(1,1),\ldots,(1,1)}_{d}),$$
we write in particular
$$\langle\tau_1^{a_1}(pt)\cdots\tau_m^{a_m}(pt)\rangle^{\text{trop}}_{d,0}:=\widetilde{N}_{\Delta,\mathbf{k}}^\text{trop}(\alpha)
=\langle\tau_1^{a_1}(pt)\cdots\tau_m^{a_m}(pt)\rangle_{d,0},$$ where the last equality is Remark \ref{rem-old_descendant}.
\end{notation}
\begin{definition}[Tropical refined descendant curve]\label{def-trdi}
Let $\Delta=(v(y_1),\ldots, v(y_n))$ and $F\subset \{1,\ldots,n\}$. An $(r,s)$-marked curve $C\in \Mrs(\Delta)$ with set $F$ of fixed ends is called \textit{(tropical) refined descendant curve} if its markings satisfy the following conditions
\begin{itemize}
 \item each of the  $r$ real markings has to be adjacent to a $3$-valent vertex of the underlying graph $\Gamma$,
 \item each of the $s$ complex markings has to be adjacent to a $4$-valent vertex of 
$\Gamma$. 
\end{itemize}
\end{definition}
\begin{example}
Note that not all curves in $\Mrs(\Delta)$ are refined descendant curves. For example, Welschinger and bridge curves containing a vertex of type (7) do not satisfy the conditions of Definition \ref{def-trdi}. This class of curves is larger than the class of broccoli curves, see Example \ref{ex-counter}, and is also larger than the class of refined broccoli curves, since for example $4$-valent vertices without an adjacent complex marking are allowed.\\
Even if we restrict ourselves to curves passing through a configuration $\calP$ of points in general position as defined in Subsection \ref{subsub-ev}, we do not obtain all curves in $\Mrs(\Delta)$ through $\calP$, see also Example \ref{ex-gegen_gen}, and this class of 
curves is still larger than the class of (refined) broccoli curves. 
\end{example}
For the following lemma we have to make sure that the number of conditions is the right one to guarantee that we obtain a finite number of curves through the point configuration.
\begin{lemma}[Bijection between unoriented refined broccoli curves and refined descendant curves]\ \\ \label{lem-equ}
Let $r,s\geq 0$ and let $\Delta=(v_1,\ldots,v_n)$ a collection of vectors in $\ZZ^2\setminus \{0\}$, and let $F\subset \{1,\ldots,n\}$ such that $r+2s+|F|=|\Delta|-1$. Furthermore, let $\calP=(P_1,\ldots,P_{2(r+s)+|F|})$  be a configuration of points in general position for $\ev_F$. Then there is a bijection between the unoriented refined broccoli curves and refined descendant curves in $\Mrs(\Delta)$ through $\calP$.
\end{lemma}

\begin{proof}
Given a configuration $\calP$ of points in general position, then clearly every unoriented refined broccoli curve also satisfies the conditions of Definition \ref{def-trdi}. Assume now that $C$ is a refined descendant curve passing through $\calP$, which is not an unoriented refined broccoli curve. We want to arrive at a contradiction.
We can assume without loss of generality that $C$ has a $4$-valent vertex to which no marking is adjacent and $F=\emptyset$. Consider first the case where $s=0$. In this situation our definition of points in general position agrees with the one of \cite{GM05b} and therefore whenever we have a $4$-valent vertex in $C$ the points of $\calP$ cannot lie in general position.\\
In order to understand the case $s \neq 0$ let us classify curves $C$ passing through $\calP$ having one $4$-valent vertex when $s=0$. Following the argumentation of \cite[Remark 3.22]{GS11}, we see that such a curve has exactly one region depicted below with dotted lines containing the $4$-valent vertex and is of type (B) or (C). By a region we understand a connected component of $\Gamma \setminus \cup_{i=1}^r \overline{x_i}$, where $\Gamma$ is the underlying graph of $C$. 
\begin {center} \input {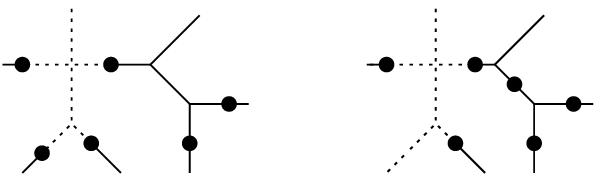} \end {center}
In case (B) all regions have exactly one end each and in case (C) all regions have one end, except that  the region with the $4$-valent vertex has two ends and there is one region having no end at all. The transition $s=0 \rightarrow s=1, r \rightarrow r-2$ modifies these graphs, namely two of the real markings on edges are replaced by a complex marking on a $3$-valent vertex. The number $b$ of bounded edges in $C$ decreases by $2$ during this transition. But already before the transition $b$ was equal to $2r-3$ in both cases (B) and (C), so following Section \ref{subsub-nonorient}, the combinatorial type has in both cases  dimension $2r-1$. After the transition, the combinatorial type has dimension $2r-3$ which is smaller than the dimension of a maximal cell in the moduli space which has dimension $2(r-2)+2s=2r-4+2=2r-2$.
\end{proof}

The previous lemma allows to make the following definition.

\begin{definition}[Tropical refined descendant invariants]\label{def-desc-inv}
Let $\Mrs^{desc}(\Delta)$ be the closure of refined descendant curves in $\Mrs(\Delta)$ and 
assume we have $r+2s+|F|=|\Delta|-1$. If $\calP$ is a collection of conditions in general position for $\ev_{\Mrs^{desc}(\Delta)}: \Mrs^{desc}(\Delta) \to \RR^{2(r+s)}$, then the \textit{tropical refined descendant invariant} is defined as $$\Nrs^{desc}(y,\Delta,F,\calP)=\frac{1}{|G(\Delta,F)|}\sum_{C}m_C(y),$$ where the sum is taken over all refined descendant curves $C$ in $\Mrs^{desc}(\Delta)$ with $\ev_{\Mrs^{desc}(\Delta)}(C)=\calP$, $m_C(y)$ is the multiplicity of the by Lemma \ref{lem-equ} associated refined unoriented broccoli curve, and $G(\Delta,F)$ is defined in Definition \ref{def-wbb}.
\end{definition}

\begin{corollary}[To Lemma \ref{lem-equ}] \label{cor-same}
Let $\calP$ be a collection of conditions in general position for $\ev_{M^{desc}_{r,s}(\Delta)}$ and for $\ev_{M^{rB}_{(r,s)}(\Delta,F)}$.
Then $$\Nrs^{desc}(y,\Delta,F,\calP)=\Nrs^{rB}(y,\Delta,F,\calP).$$
\end{corollary}
The following definition of points in general position is closer to \cite{GM05b} and \cite{MR09} than the one given in Subsection \ref{subsub-ev}. It will turn out to be useful for the last lemma of this section. In fact, our definition of tropical descendant invariants also allows  vertices in curves $C$ which have valence bigger than $3$, that is, we do not consider only zeroth and first order descendant Gromov-Witten invariants. But for refined broccoli invariants such curves are not allowed! 
\begin{definition}[Restricted general position of collections of conditions]\label{def-restricted-pts}
Let $C$ be an $(r,s)$-marked curve of degree $\Delta=(v(y_1),\ldots,v(y_n))$ and set of fixed ends $F$, and let $M\subset \Mrs(\Delta)$ be a polyhedral subcomplex. A collection $\calP=\big((P_1,\ldots,P_{r+s}),(Q_i:~i\in F)\big)$ of conditions for $\ev_{F,M}$ is called \textit{conditions in restricted general position for} $\ev_{F,M}$ if it is in general position in the sense of \ref{subsub-ev}, but where we allow only combinatorial types $\alpha$ such that every vertex in the underlying graph $\Gamma$ of a  (and therefore every) curve $C$ of this type $\alpha$ in $\ev_{F,M}^{-1}(\calP)$, that is not $3$-valent, has to be $4$-valent with a complex marking adjacent to it.
\end{definition}
\begin{lemma}[Relationship of tropical refined descendant to tropical descendant Gromov-Witten invariants]\ \\
\label{dGW=rdGW}
Under the assumptions and with the notations of Definition \ref{def-desc-inv}, let as in Definition \ref{def-dGW} $\alpha=(\alpha_1,\alpha_2,\ldots)$ be the sequence defined by 
$\alpha_i=\#\big\{v(y_j)\in \Delta\bigm | j\in F \hbox{ and $v(y_j)$ has weight }i\big\}$.
Then 
$$N^{desc}_{(r,s)}(1,\Delta,F,\calP)=\frac{\prod_{i\in O_f} w(y_i)}{\prod_{i\in O_n} w(y_i)}  I^\alpha \widetilde N^{\text{trop}}_{\Delta,(r,s,0,0,\ldots)}(\alpha),$$
where $O_f=\big\{i\in F\bigm| w(y_i)\hbox{ odd}\big\}$ and $O_n=\big\{i\not \in F\bigm| w(y_i)\hbox{ odd}\big\}$.
\end{lemma}
\ignore{\marginpar{NOT ANYMORE: what troubles me a bit is that for the Psi-count the dimension of a marking at 
a $4$-valent vertex is $1$, while in the broccoli setting each complex marking 
contributes with $2$ in the dimension count.-> WRITE SOMETHING about it!}}
\begin{proof}
Using Lemma \ref{lem-equ} we know that refined descendant curves through the right number of points in general position are composed of vertices (I')-(III'). Expanding the refined vertex multiplicity of a vertex of type (II') as in  the proof of Lemma \ref{lem-welldefined} and plugging in $y=1$ gives $a$, and plugging in $y=1$ in the refined multiplicity of a vertex of type (III') gives directly $1$.\\ 
Now, note that in the definition of tropical descendant Gromov-Witten invariants there is a factor of $\frac{1}{I^\alpha}$ introduced to compare them more easily with classical descendant Gromov-Witten invariants, whereas refined descendant invariants do not come with this factor. The factor $\frac{\prod_{i\in O_f} w(y_i)}{\prod_{i\in O_n} w(y_i)}$ should be included in the equality above since tropical refined descendant invariants are defined via refined broccoli invariants and the multiplicity of ends $y_i$ of weight $w(y_i)>1$ as defined in Definition \ref{def-or} is not always equal to $1$, when specializing $y=1$. Indeed we have for $y_i \in F$ and at $y=1$:
$$m_{y_i}(y) = \begin{cases} w(y_i) &\mbox{if } w(y_i) \text{ is odd} \\
1 & \mbox{if } w(y_i) \text{ is even,} \end{cases}$$
and for $y_i \notin F$ and at $y=1$:
$$m_{y_i}(y) = \begin{cases} \frac{1}{w(y_i)} &\mbox{if } w(y_i) \text{ is odd} \\
1 & \mbox{if } w(y_i) \text{ is even.} \end{cases}$$
\end{proof}
\begin{remark}
Geometrically speaking what happens during the transition from refined descendant broccoli curves to curves counting for tropical descendant Gromov-Witten invariants is that we replace each (big) complex marking adjacent to $4$-valent vertex of $\Gamma$ by a (small) real marking since for descendant curves there is no distinction between small and big markings, that is we have $m=r+s$ in Definition \ref{def-dGW}.
\end{remark}

\begin{remark}\label{otherdes}
We can define an alternative version $N^{desc,*}_{(r,s)}(y,\Delta,F,\calP)$ of refined descendent invariants by replacing the multiplicities $m_C(y)$ by the simpler multiplicities $m'_C(y)$ from Remark \ref{other}. Then we get $N^{desc,*}_{(r,s)}(1,\Delta,F,\calP)=  I^\alpha \widetilde N^{\text{trop}}_{\Delta,(r,s,0,0,\ldots)}(\alpha).$
\end{remark}
  \section {Proof of invariance for refined broccoli numbers} \label {sec-inv}
In this section we will show the invariance of the oriented refined broccoli invariant $\Nrs^{rB}(y,\Delta,F,\calP)$ as defined in Definition \ref{def-or-refined-brocc-inv}, i.e. we show that it does not depend on the choice of the point configuration $\calP$ as long as points in $\calP$ are in general position. The idea of proof is the same as in the proof of \cite[Theorem 3.6]{GMS13}.\\ 
Namely, we have to prove that the function $\calP \mapsto \Nrs^{rB}(y,\Delta,F,\calP)$ is locally constant on the open subset of $\RR^{2(r+s)+|F|}$ of conditions in general position for refined oriented broccoli curves, and can only differ at the image under $\ev_F$ of the boundary of top-dimensional cells of the closure $\Mrs^{rB}(\Delta,F)$ of oriented refined broccoli curves in $\Mrs^\ori(\Delta,F)$. It turns out that it is enough to show that the function $\calP \mapsto \Nrs^{rB}(y,\Delta,F,\calP)$ is locally constant around a cell in this image of codimension $1$ in $\RR^{2(r+s)+|F|}$, since any two top-dimensional cells of $\RR^{2(r+s)+|F|}$ can be connected to each other through codimension-1 cells. If $\alpha$ is a combinatorial type in $\Mrs^{rB}(\Delta,F)$ of codimension 1 such that $\ev_F$ is injective on $\Mrs^{\alpha}(\Delta,F)$, then $\ev_F$ is piecewise linear, hence it maps this cell to a unique hyperplane $H$ in $\RR^{2(r+s)+|F|}$. If furthermore $U_\alpha \subset \Mrs^{rB}(\Delta,F)$ is the open subset consisting of $\Mrs^\alpha(\Delta,F)$ together with all top-dimensional adjacent cells of $\Mrs^{rB}(\Delta,F)$, then we have to prove for a point configuration $\calP$ in a neighborhood of $\ev_F(\Mrs^\alpha(\Delta,F))$ that the sum of the multiplicities of the curves in $U_\alpha \cap \ev_F^{-1}(\calP)$ does not depend on $\calP$, that is, it is the same on both sides of $H$. More precisely, we take a configuration $\calP_+$ on one side of the hyperplane $H$ and a configuration $\calP_-$
on the other side. Then, we have to show \begin{equation}\label{wallcross} \sum_{C\in A_+} m_C(y)=\sum_{C\in A_-} m_C(y),\end{equation} where $A_{\pm}$ denotes the curves through $\calP_{\pm}$. We can take care of the side of $H$ on which a curve $C$ lies by assigning to the curve a sign $\sigma_C \in \{\pm 1\}$, called the \textit{$H$-sign}. More detail about the $H$-sign can be found in the proof of \cite[Theorem 3.6]{GMS13}. The equation (\ref{wallcross}) is equivalent to 
\begin{equation}\label{wallcrossH} \sum_{C} \sigma_C m_C(y)=0,\end{equation}
where now $C$ runs through all curves through $\calP_+$ and all curves through $\calP_-$.

\begin{theorem}[Invariance of oriented refined broccoli invariants]\label{thm-invariance}
Let $\Delta \subset \RR^2$ be a lattice polygon, let $\Mrs^{rB}(\Delta,F)$ be the closure of oriented refined broccoli curves in $\Mrs^{\text{or}}(\Delta,F)$ and assume we have $r+2s+|F|=|\Delta|-1$. Furthermore let $\calP$ be a collection of conditions in general position for $\ev_{\Mrs^{rB}(\Delta,F)}: \Mrs^{rB}(\Delta,F) \to \RR^{2(r+s)+|F|}$. Then the refined broccoli invariant $\Nrs^{rB}(y,\Delta,F,\calP)$ does not depend on $\calP$, i.e.\ it is indeed an invariant.\\ 

Therefore, we write $\Nrs^{rB}(y,\Delta,F)$ instead of $\Nrs^{rB}(y,\Delta,F,\calP)$ in the following.
\end{theorem}

\begin{proof}
Below, one can find a list of codimension-$1$ combinatorial types $\alpha$ in $\Mrs^{rB}(\Delta,F)$. In these codimension-1 cases the parity of the edges adjacent to the vertices does not play any role since contrary to the situation in \cite{GMS13} for old broccoli curves the multiplicity of a vertex in a refined oriented broccoli curve does not depend on the parity of the edges adjacent to it. Note that (in this spirit) in the pictures below, the parity of the edges can be any as long as the balancing condition is satisfied.\\
The labeling of the edges is chosen as follows. In the cases (A) and (C) there is only one edge with an orientation towards the vertex which we will denote by $v_1$. The other two, respectively three edges are labeled clockwise as $v_2,~v_3$, respectively $v_2,~v_3,~v_4$. Analogously, in case (B) there is exactly one edge with an orientation away from the vertex which we name $v_1$ and the other three are labeled clockwise as $v_2,~v_3,~v_4$. Note that the choice of resolutions below depends on the labeling and also the form of the relation below in each case depends on the labeling! Also, note that the labeling of the edges in each case is the same as in \cite{GMS13}. 
\begin{itemize}
 \item[(A)] a vertex (I) is merging with a vertex (II), which leads to a $4$-valent vertex with one real marking, two outgoing edges, and one incoming edge.
\item[(B)] a vertex (II) is merging with another vertex of type (II), yielding a $4$-valent vertex with no marking, one outgoing edge, and three incoming edges.
\item[(C)] a vertex (II) is merging with a vertex of type (III), giving rise to a $5$-valent vertex with one complex marking, three outgoing edges, and one incoming edge.
\end{itemize}
\begin{center}
\begin{picture}(0,0)%
\includegraphics{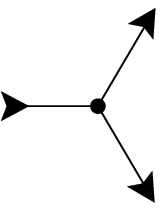}%
\end{picture}%
\setlength{\unitlength}{4144sp}%
\begingroup\makeatletter\ifx\SetFigFont\undefined%
\gdef\SetFigFont#1#2#3#4#5{%
  \reset@font\fontsize{#1}{#2pt}%
  \fontfamily{#3}\fontseries{#4}\fontshape{#5}%
  \selectfont}%
\fi\endgroup%
\begin{picture}(723,1101)(588,-2051)
\put(991,-1996){\makebox(0,0)[b]{\smash{{\SetFigFont{10}{12.0}{\familydefault}{\mddefault}{\updefault}{\color[rgb]{0,0,0}(A)}%
}}}}
\end{picture}%
\hspace{2cm}\begin{picture}(0,0)%
\includegraphics{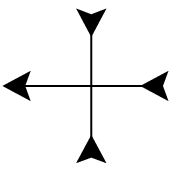}%
\end{picture}%
\setlength{\unitlength}{4144sp}%
\begingroup\makeatletter\ifx\SetFigFont\undefined%
\gdef\SetFigFont#1#2#3#4#5{%
  \reset@font\fontsize{#1}{#2pt}%
  \fontfamily{#3}\fontseries{#4}\fontshape{#5}%
  \selectfont}%
\fi\endgroup%
\begin{picture}(774,998)(484,-3536)
\put(856,-3481){\makebox(0,0)[b]{\smash{{\SetFigFont{10}{12.0}{\familydefault}{\mddefault}{\updefault}{\color[rgb]{0,0,0}(B)}%
}}}}
\end{picture}%
\hspace{2cm} \begin{picture}(0,0)%
\includegraphics{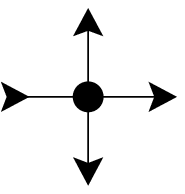}%
\end{picture}%
\setlength{\unitlength}{4144sp}%
\begingroup\makeatletter\ifx\SetFigFont\undefined%
\gdef\SetFigFont#1#2#3#4#5{%
  \reset@font\fontsize{#1}{#2pt}%
  \fontfamily{#3}\fontseries{#4}\fontshape{#5}%
  \selectfont}%
\fi\endgroup%
\begin{picture}(820,1057)(498,-4931)
\put(901,-4876){\makebox(0,0)[b]{\smash{{\SetFigFont{10}{12.0}{\familydefault}{\mddefault}{\updefault}{\color[rgb]{0,0,0}(C)}%
}}}}
\end{picture}%

\end{center}
The adjacent codimension-0 combinatorial types in $\Mrs^{rB}(\Delta,F)$ of cases $(A)$, $(B)$ and $(C)$, called \textit{resolutions} as in \cite{GMS13}, are depicted below. The vectors $v_1,\ldots,v_4$  in the figures are the outwards pointing direction vectors of the edges satisfying $v_1+v_2+v_3=0$ in case (A), and $v_1+v_2+v_3+v_4=0$  in the cases (B) and (C). Note that this means in particular that we haven chosen a labeling of the edges. Vertices in the resolutions, which are not contained in the list of Definition \ref{def-or}, are marked with a square. Resolutions containing such vertices do not contribute.
\begin {center} \input {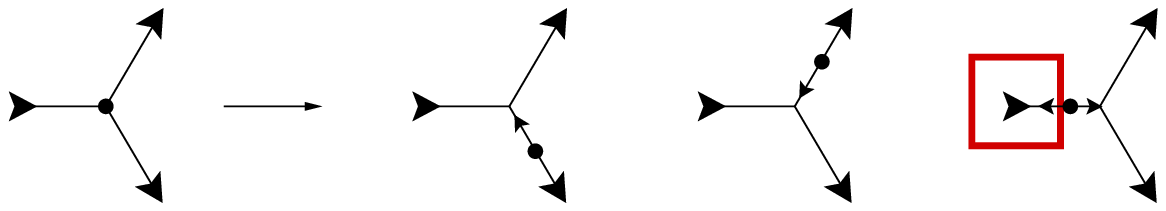} \end {center}
  \vspace {-2mm}
  \begin {center} \input {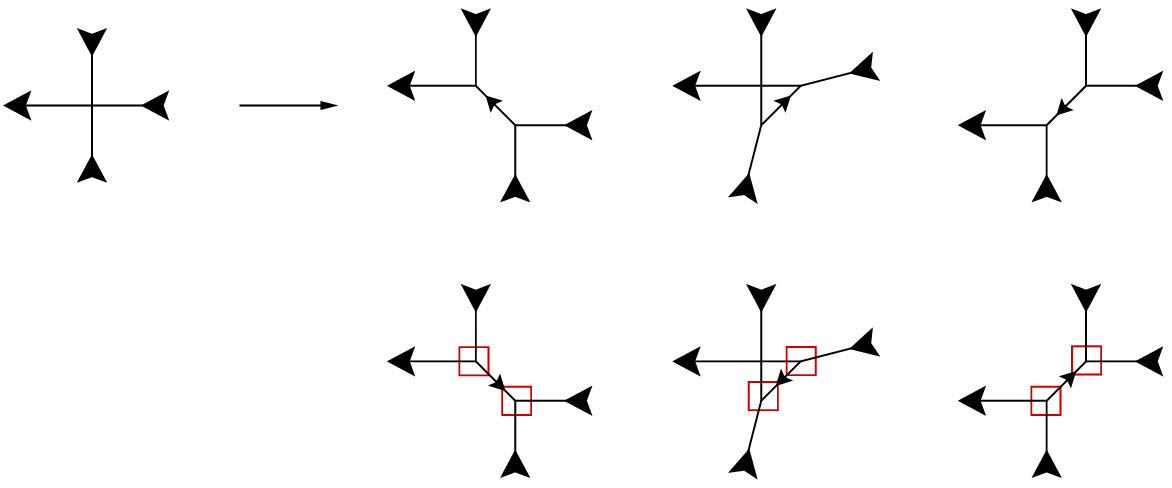} \end {center}
  \vspace {-2mm}
  \begin {center} \input {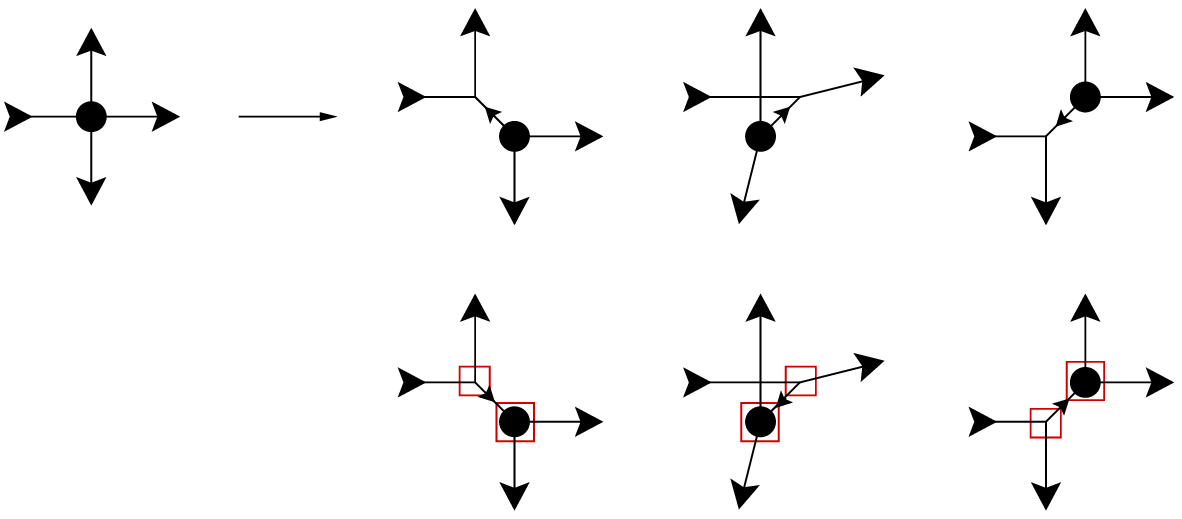} \end {center}
It is clear that there are not more possible resolutions in each case than those listed above. More precisely, in case (A) the only way to resolve the $4$-valent vertex is to put the real marking on one of the three edges of direction vector $v_i$ with $i \in \{1,2,3\}$. In cases (B) and (C) resolutions are obtained by choosing two out of four edges which first come together and then to choose the orientation of the newly created bounded edge, i.e. there are exactly 6 resolutions in each case as listed above.\\
Note also that these are exactly the same cases as in the proof of \cite[Theorem 3.6]{GMS13}. The only difference is that we do not take care of the parity of the edges in our proof and there will be in general more resolutions contributing with a multiplicity different from $0$ since more vertices are permitted. We have to prove in each case that 
\begin{equation}\label{Hres}\sum_{C}\sigma_C \cdot m_C(y)=0,
\end{equation} 
where the sum goes over all resolutions $C$ and $\sigma_C$ is the H-sign of the curve $C$ as explained in the beginning of this section. Let $(a,b)$ denote the determinant of the matrix containing the vectors $a,b \in \RR^2$ as column vectors. Following the argument of the proof of \cite[Theorem 3.6]{GMS13} we can take care of the H-sign by replacing the Mikhalkin multiplicities $a_{ij}=|(v_i,v_j)|=|(v_j,v_i)|$ of the $3$-valent vertices $V$ and $W$ in the resolutions $C$ depicted above appearing in $m_C(y)$ by the corresponding determinants $A_{ij}=(v_i,v_j)=-(v_j,v_i)$, where $(i,j)\in \{(1,2),(1,3),(1,4),(3,4),(4,2),(2,3)\}$ is the unique pair such that the $v_i$ and $v_j$  edges are adjacent to the considered vertex. 
That is \eqref{Hres} is equivalent to the formula 
\begin{equation}\label{Hres0}\sum_{C} \widetilde{m}_{C}(y)=0,
\end{equation}
where $\widetilde m_C(y)$ is obtained from $m_C(y)$ by replacing the $a_{ij}$ by the $A_{ij}$ as mentioned above.

Thus the invariance is in each case equivalent to the following equation below 
\begin{itemize}
\item[(A)] $\frac{y^{A_{12}/2}-y^{-A_{12}/2}}{y^{1/2}-y^{-1/2}}+\frac{y^{A_{13}/2}-y^{-A_{13}/2}}{y^{1/2}-y^{-1/2}}=0,$ which is true since $A_{12}+A_{13}=0$ following from the balancing condition $v_1+v_2+v_3=0$.

\item[(B)] \begin{align*}
\frac{y^{A_{12}/2}-y^{-A_{12}/2}}{y^{1/2}-y^{-1/2}}&\frac{y^{A_{34}/2}-y^{-A_{34}/2}}{y^{1/2}-y^{-1/2}}+\frac{y^{A_{23}/2}-y^{-A_{23}/2}}{y^{1/2}-y^{-1/2}}\frac{y^{A_{14}/2}-y^{-A_{14}/2}}{y^{1/2}-y^{-1/2}}\\&+\frac{y^{A_{13}/2}-y^{-A_{13}/2}}{y^{1/2}-y^{-1/2}}\frac{y^{A_{42}/2}-y^{-A_{42}/2}}{y^{1/2}-y^{-1/2}}=0.
\end{align*} Both the vertices $V$ and $W$ in all three resolutions are of the type $(II)$ of Definition \ref{def-or}, therefore they are counted with a ``-''-sign in the numerator and denominator. This equation is correct since $A_{12}+A_{13}+A_{14}=0$, $A_{12}-A_{23}+A_{42}=0$ and $A_{13}+A_{23}-A_{34}=0$, which follows from the balancing condition $v_1+v_2+v_3+v_4=0$.

\item[(C)] \begin{align*}
\frac{y^{A_{12}/2}-y^{-A_{12}/2}}{y^{1/2}-y^{-1/2}}&\frac{y^{A_{34}/2}+y^{-A_{34}/2}}{y^{1/2}+y^{-1/2}}+\frac{y^{A_{23}/2}+y^{-A_{23}/2}}{y^{1/2}+y^{-1/2}}\frac{y^{A_{14}/2}-y^{-A_{14}/2}}{y^{1/2}-y^{-1/2}}\\&+\frac{y^{A_{13}/2}-y^{-A_{13}/2}}{y^{1/2}-y^{-1/2}}\frac{y^{A_{42}/2}+y^{-A_{42}/2}}{y^{1/2}+y^{-1/2}}=0.\end{align*}
In all three resolutions  the vertex $V$ is of the type  $(II)$ of Definition \ref{def-or} and the vertex $W$ of the type $(III)$ of Definition \ref{def-or}. Therefore
$V$ is counted with $-$ sign in numerator and denominator and $W$  with a ``+''-sign in the numerator and denominator. This equation is true since $A_{12}+A_{13}+A_{14}=0$, $A_{12}-A_{23}+A_{42}=0$ and $A_{13}+A_{23}-A_{34}=0$,  which follows from the balancing condition $v_1+v_2+v_3+v_4=0$.
\end{itemize}
\ignore{
\begin{LG}
The aim is to show that this is independent of the labeling.
I am not quite sure I understand what this means:

What we see is the following:
In cases (A) and (C) the relation is clearly not independent of the labeling. The label of the incoming edge is special, and the formula depends on what the label of the incoming edge is. But the equation
 (A) is independent of the reordering of $v_2$, $v_3$ as long is $v_1$ is the incoming edge.
 (C) is independent of the reordering of $v_2,v_3,v_4$ as long as the outgoing edge is called $v_1$.
  
  So the easiest solution would be to require in (A) and (C) that the labeling will call the incoming  edge $v_1$. 
  
  Alternatively  I can show that if the incoming  edge is called $v_i$ with $i\in \{1,2,3,4\}$, then 
 with this labelling 
 (A) becomes
 
 $$\sum_{j\in \{1,2,3\}\setminus  i} \frac{y^{A_{ij}/2}- y^{-A_{ij}/2}}{y^{1/2}-y^{-1/2}}.$$
 
For $(C)$ the point is that the $-$ term is always the one containing the incoming edge.
It can be precisely determined what happens but I am not quite sure what we want. It is also a bit more complicated.

In case (B) the left hand side of the equation depends via a sign on the labelling (as the relation is that this is equal to zero, it is the same relation).
In fact if we apply a permutation to $v_1,v_2,v_3,v_4$, then both sides of the equation are just multiplied by the sign of the permutation. 

Let me explain for (B). Let $w_1,w_2,w_3,w_4$ be a permutation of $v_1,v_2,v_3,v_4$. 
Write $B_{ij}=(w_i,w_j)$
Then the argument of the paper gives the same relation $(B)$ with each $A_{ij}$ replaced by $B_{ij}$.
If the permutation is a transpositon, one sees directly that this relation, when expressed in the $A_{ij}$ is $-$ the 
original $(B)$ relation. It follows that for any permutation both sides of the relation are multiplied by 
the sign of the permutation.
\end{LG}
}

\end{proof}

\begin{remark}[Refinement of bridges]
While trying to refine the multiplicities of bridge curves as defined in Definition 5.2 \cite{GMS13}, problems occured while trying to prove the invariance along bridges similar to the proof of Theorem 5.14 \cite{GMS13}. The refined version of the bridge vertex (9) in the list \ref{def-vertextypes} as depicted below
\begin {center} \begin{picture}(0,0)%
\includegraphics{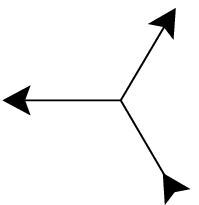}%
\end{picture}%
\setlength{\unitlength}{4144sp}%
\begingroup\makeatletter\ifx\SetFigFont\undefined%
\gdef\SetFigFont#1#2#3#4#5{%
  \reset@font\fontsize{#1}{#2pt}%
  \fontfamily{#3}\fontseries{#4}\fontshape{#5}%
  \selectfont}%
\fi\endgroup%
\begin{picture}(871,913)(2149,-630)
\end{picture}%
 \end {center}
should then have the multiplicity $$\frac{y^{a/2}+y^{-a/2}}{y^{1/2}+y^{-1/2}},$$ where $a$ is the Mikhalkin multiplicity of the vertex. 

\ignore{In case (D1) the corresponding version of the relation of Theorem \ref{thm-invariance} would be
\begin{align*}
  \frac{y^{A_{12}/2}+y^{-A_{12}/2}}{y^{1/2}+y^{-1/2}}&\frac{y^{A_{34}/2}-y^{-A_{34}/2}}{y^{1/2}-y^{-1/2}}+\frac{y^{A_{13}/2}+y^{-A_{13}/2}}{y^{1/2}+y^{-1/2}}\frac{y^{A_{24}/2}-y^{-A_{24}/2}}{y^{1/2}-y^{-1/2}}+\\&\frac{y^{A_{23}/2}+y^{-A_{23}/2}}{y^{1/2}+y^{-1/2}}\frac{y^{A_{14}/2}-y^{-A_{14}/2}}{y^{1/2}-y^{-1/2}}=0.
\end{align*}
}
Let us consider the codimension-1 case (D1) appearing in Theorem 5.14 \cite{GMS13}. As we are interested in a refined version of the bridge algorithm we should in particular choose those resolutions which have been already considered in the proof of Theorem 5.14 \cite{GMS13}, that is, resolutions IV, II and III in the picture below. Note that the sum of their refined multiplicities is $0$.
\begin {center} \input {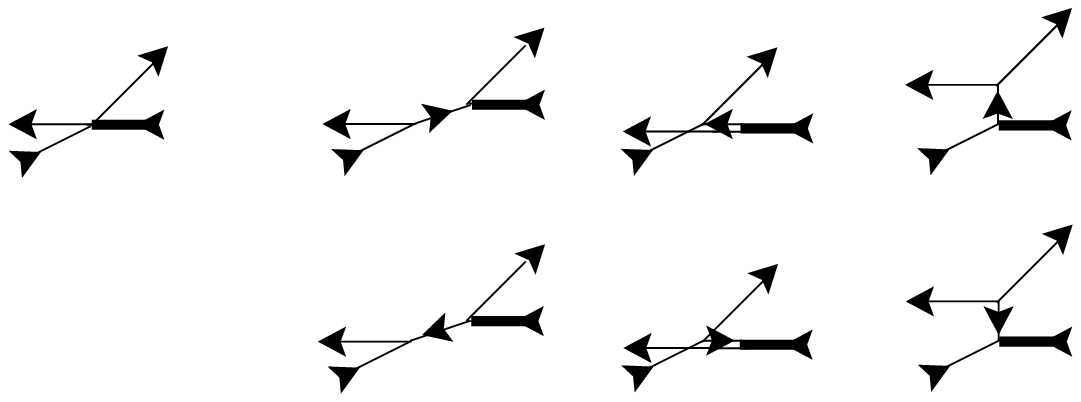} \end {center}
With the generalization of the bridge vertex (9) it can happen that the movement of a string in a curve on a bridge is not bounded. This was  excluded in Lemma 5.15 \cite{GMS13} for the standard bridge curves, if $\Delta$ is a toric Del Pezzo degree, i.e.\ the toric surface $X(\Delta)$ is $\PP^2,\PP^1 \times \PP^1$ or $\PP^2$ blown up in up to three points. For instance in the picture below we consider as toric surface $\PP^2$ and tropical curves of degree $3$. The curve below is an honest bridge curve in the sense of \cite{GMS13} and contains the two vertices from resolution III depicted above and a string (dashed) which can be moved up without hitting any other vertex of the curve. The reason for this to be possible is that the curve contains a vertex of the generalized type (9) (in the  box) with only odd edges such that the argument of Lemma 5.15 \cite{GMS13} does not work anymore.
\begin {center} \begin{picture}(0,0)%
\includegraphics{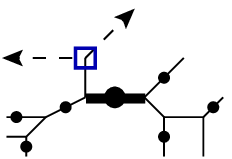}%
\end{picture}%
\setlength{\unitlength}{4144sp}%
\begingroup\makeatletter\ifx\SetFigFont\undefined%
\gdef\SetFigFont#1#2#3#4#5{%
  \reset@font\fontsize{#1}{#2pt}%
  \fontfamily{#3}\fontseries{#4}\fontshape{#5}%
  \selectfont}%
\fi\endgroup%
\begin{picture}(1033,699)(4650,-5383)
\end{picture}%
 \end {center}
This counterexample shows that, if we allow the generalized vertex type (9), the bridge algorithm does not terminate. One could hope to obtain a refined version of the bridge algorithm when one does not generalize type (9). But this does not work. Namely, we have to choose resolutions IV, II and III as in the proof of Theorem 5.14 \cite{GMS13}. It is clear that we want to choose at most three resolutions since the resolution I differs from resolution IV only by the choice of orientation of the bounded edge (analogously for the pairs II/V and III/VI), but the orientation of that edge is uniquely determined by the bridge algorithm. Now, if we do not allow the generalized vertex type (9) then resolution III counts $0$, whereas resolutions IV and II contain the original vertex type (9), respectively. Unfortunately, considering solely the resolutions IV and II does not give invariance when considering refined multiplicities.    
\end{remark}
\ignore{
\begin{remark}[Refined Welschinger invariant]\label{cor-inv_welsch}
Let $\Delta$ be a toric Del Pezzo degree as in Definition 4.22 \cite{GMS13}, which implies that the toric surface $X(\Delta)$ associated $\Delta$ equals $\PP^2$, $\PP^1\times \PP^1$ or $\PP^2$ blown up in up to three points. Let $\calP$ be a collection of conditions in general position for $\ev_{\Mrs^{rB}(\Delta,F)}: \Mrs^{rB}(\Delta,F) \to \RR^{2(r+s)+|F|}$. In this situation we can define the \textit{refined Welschinger number} as $N^{rW}_{(r,s)}(\Delta,F,\calP):=N^{rB}_{(r,s)}(\Delta,F,\calP)$, which does not depend on $\calP$ by Theorem \ref{thm-invariance}. This definition is motivated by Corollary 5.16 \cite{GMS13} where it is noted that the Welschinger numbers $N^W_{(r,s))}(\Delta,F,\calP)$ equals the corresponding broccoli number $N^B_{(r,s)}(\Delta,F,\calP)$ for a toric Del Pezzo degree $\Delta$. Nevertheless, in our situation this is only a definition and not a consequence of proven facts since we do not have a refined version for multiplicities of a bridge or Welschinger curve \marginpar{haben wir das irgendwo ordentlich definiert?}yet. 
\end{remark}
}

\begin{remark}[Refinement of multiplicities of Welschinger curves]
Although there is no obvious way to refine bridge curves such that the corresponding refined bridge algorithm terminates, it is possible to refine Welschinger curves. The refined multiplicity of a Welschinger curve then has the property that it specializes to its non-refined multiplicity as considered in \cite{GMS13}, but we cannot prove that the refined Welschinger number defined similar to the (oriented) refined broccoli number does not depend on the chosen general points in the point configuration (for this it seems that we need a bridge algorithm since also refined Welschinger numbers are not locally invariant). Remember that compared to broccoli curves Welschinger curves can also contain a vertex of type (8). If we want to refine the multiplicity associated to a Welschinger curve in the sense of Definition 4.6 \cite{GMS13} we should find in particular a refined version for the multiplicity of a vertex of type (8). It does not make sense to assign to this vertex the same refined version as to a vertex of type (4) because then the specialization to $y=-1$ will be automatically wrong. Consider the following Welschinger curve which has been already considered in the introduction of \cite{GMS13}.
\begin {center} \begin{picture}(0,0)%
\includegraphics{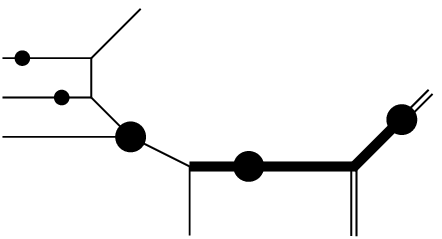}%
\end{picture}%
\setlength{\unitlength}{4144sp}%
\begingroup\makeatletter\ifx\SetFigFont\undefined%
\gdef\SetFigFont#1#2#3#4#5{%
  \reset@font\fontsize{#1}{#2pt}%
  \fontfamily{#3}\fontseries{#4}\fontshape{#5}%
  \selectfont}%
\fi\endgroup%
\begin{picture}(1989,1063)(79,-662)
\end{picture}%
 \end {center}
This curve has two vertices of type (1) contributing each $1$ to the curve multiplicity, two vertices of type (2) contributing each $1$, one vertex of type (3) contributing $2i$, one vertex of type (5) contributing $1$, one vertex of type (6b) contributing $i^{-1}$, one vertex of type (7) contributing $1$, and one vertex of type (8) contributing $-4$. So in total this curve has numerical multiplicity $-8$. Replacing the multiplicity of the vertex of type (8) by $\frac{y^2-y^{-2}}{y^{1/2}-y^{-1/2}}\cdot \frac{2}{y^{1/2}+y^{-1/2}}$ (so we consider this vertex as a vertex of type (4) times the factor $\frac{2}{y^{1/2}+y^{-1/2}}$ for the pair of edges of weight $1$) we obtain as curve multiplicity $4(y+y^{-1})$ which is equal to $-8$ when specializing to $y=-1$.\\
The easiest definition of a \textit{refined Welschinger curve} could be that is it composed of the same vertex types as non-refined Welschinger curves as defined in Definition \ref {def-wbb}, where a vertex of type 
\begin{itemize}
\item (1) counts with multiplicity $1$,
\item (2), (3) or (4) is a special case of a vertex of type (II) as defined in Definition \ref{def-or} and counts therefore with multiplicity $\frac{y^{a/2}-y^{-a/2}}{y^{1/2}-y^{-1/2}}$,
\item (5) or (6b) is a special case of a vertex of type (III) as defined in Definition \ref{def-or} and counts therefore with multiplicity $\frac{y^{a/2}+y^{-a/2}}{y^{1/2}+y^{-1/2}}$,
\item (7) counts with multiplicity $1$,
\item (8) counts with multiplicity $\frac{2(y^{a/2}-y^{-a/2})}{y-y^{-1}}$.
\end{itemize}
It would be interesting to know whether with these definitions, the count of refined Welschinger curves is indeed independent of the point conditions. In this case these "refined Welschinger invariants" should specialize to the Welschinger invariants at $y=-1$, and thus equal the specialization of the refined broccoli invariants for toric del Pezzo degrees \cite[Corollary 5.16]{GMS13}. 
As the refined bridge algorithm does not terminate, it is not clear whether in this case the refined Welschinger invariants would be equal to the refined broccoli invariants, leaving room for the more exciting possibility, that in this way new refined invariants can be obtained, maybe with a different enumerative meaning at $y=1$. Nevertheless, it is clear that we do not have local invariance for refined Welschinger invariants using the definition of a refined Welschinger curve from above, which can be seen easily from the following example. Consider the resolution of the codimension-$1$ case depicted below.
\begin {center} \begin{picture}(0,0)%
\includegraphics{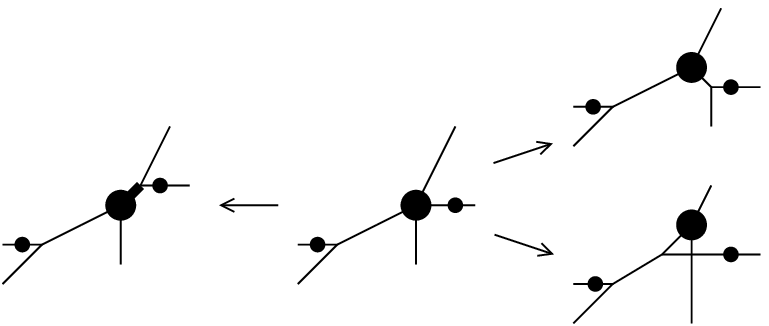}%
\end{picture}%
\setlength{\unitlength}{4144sp}%
\begingroup\makeatletter\ifx\SetFigFont\undefined%
\gdef\SetFigFont#1#2#3#4#5{%
  \reset@font\fontsize{#1}{#2pt}%
  \fontfamily{#3}\fontseries{#4}\fontshape{#5}%
  \selectfont}%
\fi\endgroup%
\begin{picture}(3489,1464)(-461,-1693)
\end{picture}%
 \end {center}
The curve on the left-hand side is a broccoli curve, but not a Welschinger curve and hence not a refined Welschinger curve. It counts with multiplicity $0$. Both curves on the right-hand side are Welschinger and also refined broccoli  curves. The curve above counts with (refined) multiplicity $y^1-1+y^{-1}$ and the curve below with (refined) multiplicity $1$. Since the sum of multiplicities $y^1+y^{-1}$ is different from $0$ the claim follows.   
\end{remark}

\begin{remark}
Consider the vertex types from Definition \ref{def-nonor}. Both for refined (unoriented) broccoli curves and for (non-refined) descendant curves a vertex of type (I') contributes with multiplicity $1$, as does for descendant
curves
a vertex of type (III'). A vertex of type (II') contributes with multiplicity $\frac{y^{a/2}-y^{-a/2}}{y^{1/2}-y^{-1/2}}$ for a refined broccoli curve, but for a descendant curve just with the Mikalkin multiplicity $a$. Similarly for a vertex of type (III'): it contributes with multiplicity $\frac{y^{a/2}+y^{-a/2}}{y^{1/2}+y^{-1/2}}$ for a refined broccoli curve, but with multiplicity $1$ for a descendant curve. We have invariance of rational descendant numbers (Theorem 8.4. \cite{MR09}) and refined broccoli numbers (Theorem \ref{thm-invariance}) respectively. Furthermore, we know by Corollary \ref{cor-same} and Lemma \ref{dGW=rdGW} that the specialization $y=1$ for a refined broccoli invariant gives the corresponding rational descendant invariant if we consider point configurations containing only points in restricted general position. Nevertheless, it is not possible to mix the vertex multiplicities for vertices in refined broccoli curves and descendant curves in order to produce an invariant. Define the following cases
\begin{itemize}
 \item[(a)] Count vertices of type (II') with Mikalkin multiplicity $a$. This is the specialization at $y=1$ of the multiplicity in (a'). 
 \item[(a')] Count vertices of type (II') with multiplicity $\frac{y^{a/2}-y^{-a/2}}{y^{1/2}-y^{-1/2}}$.
\item[(b)] Count vertices of type (III') with multiplicity $1$. This is the specialization $y=1$ of the multiplicity in (b'). 
\item[(b')] Count vertices of type (III') with multiplicity $\frac{y^{a/2}+y^{-a/2}}{y^{1/2}+y^{-1/2}}$.
\end{itemize}
For the refined broccoli invariants we use the combination $(a')+(b')$, it specializes at $y=1$ to $(a)+(b)$, which is the combination of the rational descendant invariants. The following example shows that the combination (a')+(b) does not give the expected curve multiplicity for the specialization $y=-1$. Consider degree $2$ curves through a point configuration of three real points and one complex point ($r=3$, $s=1$). If the points are ordered as depicted on the left hand side of the figure below, there is exactly one curve through the point configuration (considering unlabeled curves), which counts with multiplicity $1$. Similarly for the curve on the right hand side, which counts $y^{1/2}+y^{-1/2}$, but as we consider here unlabeled curves, we have to include an automorphism factor of $1/2$. Obviously, these two multiplicities do not agree. As there is each time exactly one curve through the configuration, this implies that we do not obtain invariance for this choice of vertex multiplicities. 
\begin {center} \begin{picture}(0,0)%
\includegraphics{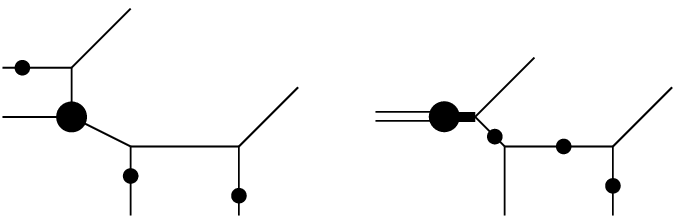}%
\end{picture}%
\setlength{\unitlength}{4144sp}%
\begingroup\makeatletter\ifx\SetFigFont\undefined%
\gdef\SetFigFont#1#2#3#4#5{%
  \reset@font\fontsize{#1}{#2pt}%
  \fontfamily{#3}\fontseries{#4}\fontshape{#5}%
  \selectfont}%
\fi\endgroup%
\begin{picture}(3084,969)(349,-658)
\end{picture}%
 \end {center}
The same example shows that we do no obtain invariance with the combination (a)+(b'). The curve on the left hand side always counts $1$, whereas the curve on the right hand side counts $\frac{4}{y^{1/2}+y^{-1/2}}$, times the automorphism factor of $1/2$. 
\ignore{The picture below shows an example where the combination (a')+(b)+(b') gives a number which is not an invariant. Again, we consider degree $2$ curves, but this time point configurations containing one real point and two complex points ($r=1$, $s=2$). Through each configuration depicted below there is exactly one curve. For the curve on the left hand side it does not matter if we count the vertex $A$ and/or the vertex $B$ as in (b) or in (b'), that is the curve always counts with multiplicity $1$. For the other curve the situation is different. If we count the vertex $C$ as in (b) and the vertex $D$ as in (b') we obtain $y^{1/2}+y^{-1/2}$ as curve multiplicity. In contrast, when we count the vertex $C$ as in (b') and the vertex $D$ as in (b) we get $(y^{1/2}+y^{-1/2})(\frac{2}{y^{1/2}+y^{-1/2}})=2$ \marginpar{there is no factor 1/2 in both examples, right?}as curve multiplicity. So it is also not possible to count some of the vertices of (III) as in (b) and some as (b').    
\begin {center} \input {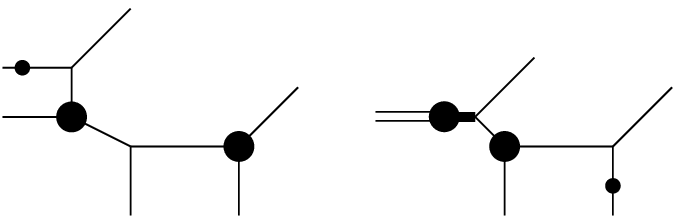} \end {center}}
\end{remark}

  \bibliographystyle {amsalpha}
  \bibliography {references}

\providecommand{\bysame}{\leavevmode\hbox to3em{\hrulefill}\thinspace}
\providecommand{\MR}{\relax\ifhmode\unskip\space\fi MR }
\providecommand{\MRhref}[2]{%
  \href{http://www.ams.org/mathscinet-getitem?mr=#1}{#2}
}
\providecommand{\href}[2]{#2}
\begin{thebibliography}{GKM09}

\bibitem[BG16]{BGo14}
F.~Block and L.~G{\"o}ttsche, \emph{{R}efined curve counting with tropical
  geometry}, Compos. Math. \textbf{152} (2016), no.~1, 115--151.

\bibitem[BGM12]{BGM11}
F.~Block, A.~Gathmann, and H.~Markwig, \emph{Psi-floor diagrams and a
  {Caporaso-Harris} type recursion}, Israel J. Math. \textbf{191} (2012),
  405--449.

\bibitem[FP97]{FP97}
William Fulton and Rahul Pandharipande, \emph{Notes on stable maps and quantum
  cohomology}, Proceedings of the 1995 Santa Cruz conference, 1997, arXiv:
  alg-geom/9608011.

\bibitem[FS15]{FS12}
S.A. Filippini and J.~Stoppa, \emph{{B}lock-{G}\"ottsche invariants from
  wall-crossing}, Compos.\ Math.\ \textbf{151} (2015), no.~8, 1543--1567.

\bibitem[GK]{GoK15}
L.~G{\"o}ttsche and B.~Kikwai, \emph{Refined node polynomials via long edge
  graphs}, to appear in Commun.~Number Theory Phys., arXiv:1511.02726.

\bibitem[GKM09]{GKM07}
A.~Gathmann, M.~Kerber, and H.~Markwig, \emph{Tropical fans and the moduli
  space of rational tropical curves}, Compos. Math. \textbf{145} (2009), no.~1,
  173--195.

\bibitem[GM07]{GM05a}
A.~Gathmann and H.~Markwig, \emph{The numbers of tropical plane curves through
  points in general position}, J. reine angew. Math. \textbf{602} (2007),
  155--177.

\bibitem[GM08]{GM05b}
\bysame, \emph{Kontsevich's formula and the {WDVV} equations in tropical
  geometry}, Adv. Math. \textbf{217} (2008), 537--560.

\bibitem[GMS13]{GMS13}
A.~Gathmann, H.~Markwig, and F.~Schroeter, \emph{Broccoli curves and the
  tropical invariance of {W}elschinger numbers}, Adv.\ Math.\ \textbf{240}
  (2013), 520--574.

\bibitem[GS]{GoSch}
L.~G{\"o}ttsche and F.~Schroeter, \emph{Floor diagrams for refined broccoli
  curves}, in preparation.

\bibitem[GS12]{GS11}
A.~Gathmann and F.~Schroeter, \emph{Irreducible cycles and points in special
  position in moduli spaces for tropical curves}, Elect.\ J.\ Comb.\
  \textbf{19} (2012), no.~4, P26.

\bibitem[GS14]{GoS12}
L.~G{\"o}ttsche and V.~Shende, \emph{Refined curve counting on complex
  surfaces}, Geom.~Topol. \textbf{18} (2014), no.~4, 2245--2307.

\bibitem[IKS09]{IKS09}
I.~Itenberg, V.~Kharlamov, and E.~Shustin, \emph{A {Caporaso-Harris} type
  formula for {Welschinger} invariants of real toric {Del Pezzo} surfaces},
  Comment. Math. Helv. \textbf{84} (2009), 87--126.

\bibitem[IM13]{IM13}
I.~Itenberg and G.~Mikhalkin, \emph{On {B}lock-{G}oettsche multiplicities for
  planar tropical curves}, Intern.\ Math.\ Res.\ Notices\ \textbf{23} (2013),
  5289 -- 5320.

\bibitem[KM98]{KM98}
M.~Kontsevich and Y.~Manin, \emph{Relations between the correlators of the
  topological sigma-model coupled to gravity}, Commun. Math. Phys. \textbf{196}
  (1998), no.~2, 385--398.

\bibitem[Mik]{Mik15}
G.~Mikhalkin, \emph{Quantum indices of real plane curves and refined
  enumerative geometry}, arXiv:1505.04338.

\bibitem[Mik05]{Mik05}
\bysame, \emph{Enumerative tropical geometry in $\mathbb{R}^2$}, J. Amer.\
  Math.\ Soc.\ \textbf{18} (2005), no.~2, 313--377.

\bibitem[Mik06]{Mik06}
\bysame, \emph{Tropical geometry and its applications}, Proceedings of the ICM
  in Madrid 2006, vol.~II, EMS, 2006, pp.~827--852.

\bibitem[MR09]{MR09}
H.~Markwig and J.~Rau, \emph{Tropical descendant {G}romov-{W}itten invariants},
  Manuscr. Math. \textbf{129} (2009), no.~3, 293--335.

\bibitem[NPS]{NPS16}
J.~Nicaise, S.~Payne, and F.~Schroeter, \emph{Tropical refined curve counting
  via motivic integration}, arXiv:1603.08424.

\bibitem[Shu06]{Shu06}
E.~Shustin, \emph{A tropical calculation of the {Welschinger} invariants of
  real toric {Del Pezzo} surfaces}, J. Algebraic Geom. \textbf{15} (2006),
  no.~2, 285--322.

\bibitem[Wel03]{Wel03}
J.-Y. Welschinger, \emph{Invariants of real rational symplectic 4-manifolds and
  lower bounds in real enumerative geometry}, C. R.\ Math.\ Acad.\ Sci.\ Paris
  \textbf{336} (2003), no.~4, 341--344.

\bibitem[Wel05]{Wel05}
\bysame, \emph{Invariants of real symplectic 4-manifolds and lower bounds in
  real enumerative geometry}, Invent.\ Math.\ \textbf{162} (2005), no.~1,
  195--234.

\end{thebibliography}

\end {document}